\documentclass[reqno]{amsart}

\usepackage{a4wide}
\usepackage{hyperref}
\usepackage{stmaryrd}
\usepackage{enumerate}
\usepackage{easybmat}
\usepackage[shortcuts]{extdash}

\renewcommand{\o}{\circ}
\def\frak{\mathfrak}
\def\Bbb{\mathbb}
\def\Cal{\mathcal}

\newcommand{\al}{\alpha}
\newcommand{\be}{\beta}
\newcommand{\ga}{\gamma}
\newcommand{\ep}{\varepsilon}
\newcommand{\eps}{\epsilon}
\newcommand{\de}{\delta}
\newcommand{\ka}{\kappa}
\newcommand{\la}{\lambda}

\newcommand{\om}{\omega}
\newcommand{\ph}{\varphi}

\newcommand{\si}{\sigma}

\newcommand{\Ga}{\Gamma}
\newcommand{\La}{\Lambda}

\newcommand{\Rho}{{\mbox{\sf P}}}
\newcommand{\Om}{\Omega}

\newcommand{\Y}{\Upsilon}
\newcommand{\ups}{\upsilon}

\newcommand{\Ad}{\operatorname{Ad}}

\renewcommand{\exp}{\operatorname{exp}}
\newcommand{\id}{\operatorname{id}}

\newcommand{\im}{\operatorname{im}}
\newcommand{\tr}{\operatorname{tr}}

\newcommand{\End}{\operatorname{End}}
\newcommand{\Aut}{\operatorname{Aut}}
\newcommand{\rank}{\operatorname{rank}}
\renewcommand{\o}{\circ}

\def\Re{\operatorname{Re}}
\def\Im{\operatorname{Im}}
\def\Gr{\operatorname{Gr}}

\def\P{{\Bbb P}}
\def\Proj{{\Cal P}}

\def\sl{\frak{sl}}

\def\gl{\frak{gl}}
\def\g{\frak g}

\def\p{\frak p}
\def\q{\frak q}

\def\r{\frak r}
\def\X{\frak X}
\def\({\big(}
\def\){\big)}
\def\R{\Bbb R}

\def\H{\Bbb H}

\def\G{{\Cal G}}

\def\Q{{\Cal Q}}
\def\Z{{\Cal Z}}
\def\J{{\Cal J}}
\def\T{{\Cal T}}
\def\omi{\om^{-1}}

\def\pmat#1{\begin{pmatrix}#1\end{pmatrix}}

\def\Pmat#1{\left(\begin{BMAT}(@){cc0c}{cc0c}#1\end{BMAT}\right)}
\def\Bmat#1{\left[\begin{BMAT}(@){cc0c}{cc0c}#1\end{BMAT}\right]}
\def\BBmat#1{\left[\begin{BMAT}(@){c0c}{c0c}#1\end{BMAT}\right]}
\def\PPmat#1{\left(\begin{BMAT}(@){c0c}{c0c}#1\end{BMAT}\right)}

\let\del=\partial
\let\[=\llbracket
\let\]=\rrbracket
\let\x=\times
\def\.{\hbox to5pt{\hss$\cdot$\hss}}
\def\span#1{\langle#1\rangle}

\let\un=\underline
\let\wh=\widehat

\def\qtext#1{\quad\text{#1}\quad}

\newtheorem{prop}{Proposition}[section]
\newtheorem{thm}[prop]{Theorem}
\newtheorem{lem}[prop]{Lemma}

\theoremstyle{definition}

\begin{document}

\title[Para-quaternionic and Grassmannian geometry]
{Interactions between para-quaternionic \\ and Grassmannian geometry}

\author{Vojt\v{e}ch \v{Z}\'adn\'\i k}
\address{Faculty of Education \\ Masaryk University \\ Po\v{r}\'\i\v{c}\'\i\ 31 \\ 60300 Brno \\ Czech Republic} 
\email{zadnik@mail.muni.cz}

\subjclass[2000]{%
53C15, %%% General geometric structures on manifolds (almost complex, contact, symplectic, almost product structures, etc.)
53A40, %%% Other special differential geometries 
53C05} %%% Connections, general theory

\keywords{%
Almost para-quaternionic structures,
almost Grassmannian structures, 
Cartan connections}

\begin{abstract}
Almost para\-/quaternionic structures on smooth manifolds of dimension $2n$ are equivalent to almost Grassmannian structures of type $(2,n)$.
We remind the equivalence and exhibit some interrelations between subjects that were previously studied independently from the para\-/quaternionic and the Grassmannian point of view.
In particular, we relate the respective normalization conditions, distinguished curves and twistor constructions.
\end{abstract}

\maketitle

%%%
\section{Introduction}  \label{Intro}
Almost para-quaternionic structures are geometric structures that are related to the algebra of para-quaternions similarly as almost quaternionic structures are related to usual quaternions.
Both these structures can be seen as different real forms of the complex quaternionic structure and as such they have a lot in common.
Both these structures can also be studied from various viewpoints, which becomes apparent especially in the para-quaternionic case.
In this paper we focus on a natural equivalence between almost para\-/quaternionic structures on $2n$-dimensional manifolds and almost Grassmannian structures of type $(2,n)$.
Our main intent is to use consistently this equivalence to compare several notions and constructions that are already known and studied in respective communities, whose relationships are, however, not clearly visible in the existing literature.
This demarcates the structure of the paper.

In section \ref{Paraquat}, respectively \ref{Grassmann}, we collect basic definitions and concepts from para-quaternionic, respectively Grassmannian geometry, that are relevant to our purposes.
This involves the description of normalization conditions (and the corresponding families of compatible affine connections), distinguished curves and twistor constructions.
In contrast to almost quaternionic structures, there are distinguished (null) directions in the tangent bundle of a para-quaternionic manifold.
This fact yields a richer discussion both for distinguished curves and twistor bundles. 
In section \ref{Grassmann} we also introduce the main instrument for what follows, namely, the canonical Cartan connection that is associated to an almost Grassmannian structure.
There is nothing really original in these two sections.

In section \ref{Inter} we start with the comparisons and interactions;
the main observations of that section are summarized in Theorem \ref{p:first}.
It in particular follows, that the two normalization conditions, and hence the two families of distinguished compatible connections,  coincide.
This allows an easy account of the relation between families of distinguished curves. 

In section \ref{Twistor} we revise the twistor spaces of an almost para-quaternionic manifold, which are distinguished by the sign $\eps\in\{-1,0,1\}$, in the framework of the associated Cartan geometry. 
In this vein we can extend the known integrability results for $\eps=\pm1$ also to the case $\eps=0$, see Theorem \ref{p:intgrblt}.
Moreover, it follows that it is the 0-twistor space which provides a link between the two twistor constructions.
In particular, it allows an instant interpretation in Grassmannian terms, see Propositions \ref{p:z0} and \ref{p:z0int}.

Most of the considerations is independent of the dimension of the base manifold.
However, specific features appear in the lowest reasonable dimension, i.e. in dimension four.
That is why we have to add some remarks to this case, see section \ref{Final}. 
In the same section we also comment the situation when the structure admits a compatible metric.

\smallskip

There is wide literature both on almost para-quaternionic and almost Grassmannian (and related) structures.
For the former structures we follow primarily Alekseevsky and Cort\'es \cite{Alekseevsky2008a} and David \cite{David2009}. 
For the latter structures, our starting reference is \cite{Bailey1991} by Bailey and Eastwood.
It is worth noticing that twistor constructions discussed in these articles represent different ways of generalizing the Penrose's twistor program for four dimensional conformal structures, cf. \cite{Penrose1977}, to higher dimensions.
A predecessor and the closest relative of the former approach concerns, of course, almost quaternionic structures, connected with the work of Salamon, see e.g. \cite{Salamon1982a}.
A large amount of relevant material is incorporated in the monograph \cite{Cap2009} by \v{C}ap and Slov\'ak in the context of Cartan, respectively parabolic, geometries.
That book, especially its fourth chapter, was the main source of inspiration for this paper.

\subsection*{Acknowledgements}
I am grateful to  Dmitri V.  Alekseevsky, Andreas \v{C}ap, Jan Slov\'ak and Josef \v{S}ilhan for many helpful conversations.
This work was supported by the  Czech Science Foundation (GA\v{C}R) under the grants GA201/08/0397 and GA17-01171S at different times.

%%%
\section{Almost para\-/quaternionic structures}   \label{Paraquat}
After a quick reminder of para\-/quaternions, we describe the almost para\-/quaternionic structures, their compatible connections and related twistor constructions. 
The basic references for this section are \cite{Alekseevsky2008a}, \cite{Alekseevsky2009} and \cite{David2009}.

\subsection{Para\-/quaternions}\label{paraquat}
The algebra of \textit{para\-/quaternions}, denoted by $\H_s$, is characterized as the unique 4-dimensional real associative algebra with {indefinite} multiplicative norm.%
\footnote{Instead of the prefix para-, various synonyms can be found in the literature; split- is probably the most often.}
Para\-/quaternions are written as $q=a+bi+cj+dk$, for $a,b,c,d\in\R$, with the defining relations 
\begin{equation*}
  i^2=j^2=1 \qtext{and} k=ij=-ji.
\end{equation*}
Consequently, $k^2=-1$, $ik=-ki=j$ and $jk=-kj=-i$.
The conjugate para\-/quaternion to $q$ is $\bar q=a-bi-cj-dk$ and the norm is given by $|q|^2=q\bar q=\bar q q$; the corresponding polar form is $\span{p,q}=\Re(p\bar q)$.
Purely imaginary para\-/quaternions are characterized by
$\bar q=-q$, hence $q^2=-q\bar q=-|q|^2$ holds, for any $q\in\Im\H_s$.
The quadruple $(1,i,j,k)$ forms a real orthonormal basis of $\H_s$, where $|i|^2=|j|^2=-1$ and $|k|^2=1$.
There are two-dimensional real isotropic subspaces in $\H_s$, hence the inner product has the split signature $(2,2)$.

The algebra of para\-/quaternions 
is isomorphic to the algebra of endomorphism of $\R^2$,
i.e.  the matrix algebra $Mat_{2\x 2}(\R)$, 
such that the norm squared corresponds to the  determinant.
The isomorphism is given by
\begin{equation*}
  a+bi+cj+dk\mapsto\pmat{a+b&c+d\\ c-d&a-b}.
  \label{eq:quatmat}
\end{equation*}
In particular, the standard basis is mapped as 
\begin{equation}   
  1\mapsto\pmat{1&0\\0&1},\quad i\mapsto\pmat{1&0\\0&-1},\quad 
  j\mapsto\pmat{0&1\\1&0},\quad k\mapsto\pmat{0&1\\-1&0}.
  \label{eq:ijk}
\end{equation}
In these terms, the purely imaginary para\-/quaternions $\Im\H_s\subset\H_s$
form the three-dimensional subspace of trace-free matrices, which is invariant under the conjugation by regular matrices.

The group of automorphisms of $\H_s$ is just the subgroup of those elements of $SO(\H_s)\cong SO(2,2)$ which acts on $\Re\H_s$ by the identity.%
\footnote{Throughout the paper we use the standard notation of favourite Lie groups; $GL$ for general linear, $PGL$ for projective linear, $SL$ for special linear, $SO$ for special orthogonal, $Spin$ for spin.
The corresponding Lie algebras will be denoted as $\frak{gl}$, $\frak{sl}$ and $\frak{so}$, respectively.}
Hence $\Aut(\H_s)$ is isomorphic to $SO_0(1,2)$, the connected component of the identity element in $SO(1,2)$.
Under the identification above, the group of unit para\-/quaternions $\{q\in\H_s:|q|^2=1\}$ is isomorphic to $SL(2,\R)$.
The conjugation by any such element, $p\mapsto qpq^{-1}=qp\bar q$, yields a surjective group homomorphism $SL(2,\R)\to\Aut(\H_s)$ whose kernel is $\{\pm1\}$.
This just recovers the two-fold covering $SL(2,\R)\to PGL(2,\R)$ or, isomorphically written, $Spin(1,2)\to SO_0(1,2)$.

\subsection{Almost para\-/quaternionic structures}  \label{para}
A \textit{para\-/quaternionic structure} on a real vector space $W$ is a 3-dimensional subspace $\Q_W$ of endomorphisms of $W$ admitting a basis $(I,J,K)$ such that 
\begin{equation}   \label{eq:IJK}
  I\o I=J\o J=\id\qtext{and} K=I\o J=-J\o I.
\end{equation}
Consequently, $K\o K=-\id$, $I\o K=-K\o I=J$, $J\o K=-K\o J=-I$ and, for any $A\in\Q_W$, the composition $A\o A$ is a multiple of the identity map.
In particular, $\Q_W$ is endowed with an inner product of signature $(1,2)$, determined by 
\begin{equation} \label{eq:A^2}
  A\o A=-|A|^2\id,
\end{equation}
and an orientation, determined by any basis of $\Q_W$.
It follows the dimension of $W$ is necessarily even.
A linear map $f:W\to W$ is an \textit{para-quaternionic homomorphism} of $\Q_W$ if there is a  linear map $\phi:\Q_W\to\Q_W$ such that 
\begin{equation}\label{eq:auto}
f(AX)=\phi(A)f(X),
\end{equation}
for all $X\in W$ and $A\in\Q_W$.

An \textit{almost para\-/quaternionic structure} on a smooth manifold $M$ of even dimension $2n\ge 4$ is given by a subbundle $\Q\subset\End(TM)$ of rank 3, which is (locally) generated by the triple $(I,J,K)$ satisfying \eqref{eq:IJK}.%
\footnote{Almost para\-/quaternionic structures appear under various names in the literature; the frequent one in older references is almost quaternionic structures of second type.
There are also alternative equivalent definitions of the structure, see e.g. \cite{Yano1973} for more information.}
The bundle $\Q$ is endowed with a bundle metric \eqref{eq:A^2} so that the typical fibre of $\Q\to M$ is the standard oriented Minkowski space.
In the terminology of  \cite{Alekseevsky2008a}, \cite{Crampin1985}, 
elements $A\in\Q$ such that $|A|^2$ is 1, 0 and $-1$ 
(i.e. endomorphisms such that $A\o A$ is $-\id$, 0 and $\id$),
are called \textit{almost complex, almost tangent} and \textit{almost para-complex} structures, respectively.
We mostly use the uniform abbreviation \textit{almost $\eps$-complex structure}, where $\eps=-|A|^2\in\{-1,0,1\}$.

Note that, for an almost tangent structure $A\in\Q$, the condition $A\o A=0$ implies that $\im A=\ker A$, which yields a distinguished distribution in $TM$ of rank $n$.
Similarly, for a para-complex structure $A\in\Q$, the $(\pm 1)$-eigenspace decomposition of $TM$ forms two complementary distributions in $TM$ of the same rank $n$.
All such subspaces form a subset of distinguished elements in the tangent bundle, which is therefore an important (although often overlooked) part of the structure.
Vectors belonging to this subset are called \textit{null}.
We return to this subject in section \ref{paragrass}.

Let us emphasize that almost para\-/quaternionic manifolds may have any even dimension. 
Only the existence of a non-degenerate compatible metric brings an additional restriction so that the dimension of the base manifold has to be a multiple of four.
See section \ref{Metric} for some details.

\subsection{Para-quaternionic connections}\label{minimal}
A compatible connection of a para-quaternionic structure is a linear connection on $TM$ which preserves the subbundle $\Q\subset\End(TM)$; any such connection is called \textit{para-quaternionic}. 
Para-quaternionic connections generally have torsion that cannot be eliminated and therefore yields an important invariant of the structure.

Following the general theory of G-structures, let $m=\dim M$, let $G_0\subset GL(m,\R)$ be the structure group in question and let $\G_0\to M$ be the corresponding reduction of the principal frame bundle to $G_0$.
Then the change of compatible connection is controlled by a $G_0$-equivariant map $\G_0\to\R^{m*}\otimes\g_0$, where $\g_0$ is the Lie algebra of $G_0$.
The corresponding change of torsions is then expressed by an equivariant map $\G_0\to\im\del$, where $\del:\R^{m*}\otimes\g_0\to\La^2\R^{m*}\otimes\R^{m}$ is the composition
\begin{equation}
\R^{m*}\otimes\g_0\to\R^{m*}\otimes\R^{m*}\otimes\R^m\to\La^2\R^{m*}\otimes\R^{m},
\label{eq:del}
\end{equation}
whose first map is given by the inclusion $\g_0\subset\gl(m,\R)\cong\R^{m*}\otimes\R^m$ and the second map by the alternation.
Thus, compatible connections having the same torsion are parametrized by equivariant maps with values in 
$\ker\del=(\R^{m*}\otimes\g_0)\cap(S^2\R^{m*}\otimes\R^{m})$, the \textit{first prolongation of $\g_0$}.
A choice of $G_0$-invariant complement to $\im\del$ in $\La^2\R^{m*}\otimes\R^{m}$ may be used for the normalization of the torsion.

This way the para-quaternionic structures are discussed in \cite{David2009}.
In analogy to the case of almost quaternionic structures \cite{Alekseevsky1996}, it follows that $\ker\del\cong\R^{m*}$, i.e. that compatible connections having the same torsion are parametrized by one-forms on $M$.
Concretely, for $\Y\in\Om^1(M)$, the difference tensor of two such connections may be written as
\begin{equation}
\Y\odot\id +(\Y\o I)\odot I +(\Y\o J)\odot J -(\Y\o K)\odot K,
\label{eq:ups}
\end{equation}
where $(I,J,K)$ is a local basis of $\Q$ satisfying \eqref{eq:IJK}.
It further follows that $\del$ is not surjective and there is a natural $G_0$-invariant complement to $\im\del$ in $\La^2\R^{m*}\otimes\R^{m}$; more details are in section \ref{connections}.
Compatible connection whose torsion takes values in that complement is called \textit{minimal}.

The \textit{torsion} of the para\-/quaternionic structure is given by the projection of the torsion of any compatible connection to the just mentioned complement.
An almost para\-/quaternionic structure is called \textit{para\-/quaternionic}, or \textit{integrable}, if it has trivial torsion, i.e. if admits a torsion-free para\-/quaternionic connection.

\subsection{$Q$-planar curves}\label{Q-planar}
Almost para-quaternionic structures belong to a broad family of structures defined by a set of endomorphisms of the tangent bundle. 
As such they allow a class of distinguished curves, the generalized planar curves in the sense of  \cite{Hrdina2006}.
In our setting, a  parametrized curve $\ga:I\to M$, $I\subseteq\R$,  is called \textit{$Q$-planar} with respect to the almost para-quaternionic structure $\Q\subset\End(M)$ and a para-quaternionic connection $\nabla$ if the covariant derivative of the tangent vector field $\dot\ga$ belongs to its para-quaternionic span, i.e. if
\begin{equation*}
\nabla_{\dot\ga}{\dot\ga}=S(\dot\ga),
\end{equation*}
where $S$ is a section of $\span{\id}\oplus\Q\subset\End(TM)$ along $\ga$.
The definition is obviously independent of the parametrization of the curve.

From \eqref{eq:ups} it follows that a curve is $Q$-planar with respect to one minimal para-quaternionic connection if and only if it is $Q$-planar with respect to all of them.
Trivially, geodesics of any such connection are $Q$-planar with respect to all others, but they need not be their geodesics.
For finer discussion we have to distinguish curves that are everywhere, respectively nowhere, tangent to the subset of null elements of the tangent bundle;
the former curves are called \textit{null}, the latter \textit{generic}.
We return to this subject in section \ref{curves}.
It will, in particular, follow that any generic $Q$-planar curve is in fact geodesic of some compatible connection.
We will also specify a distinguished subclass among the class of generic $Q$-planar curves.
The discussion for null $Q$-planar curves is more strict.

\subsection{Twistor spaces for almost para\-/quaternionic manifolds}
\label{paratwist}
Given an almost para\-/quaternionic manifold $(M,\Q)$ and an arbitrary
$s\in\R$, the \textit{$s$-twistor space} $\Z^s\to M$ is defined as
$$
\Z^s:=\{A\in\Q:|A|^2=-s, \text{\ i.e.\ } A\o A=s\id\}.
$$
By definition, each $s$-twistor space is a fibre bundle over $M$ with 2-dimensional fibre, so the dimension of the total space is also even.

Following the observations of section \ref{para},
the typical fibre of $\Q\to M$ is decomposed into
disjoint subsets consisting of space-, light- and time-like vectors.
Accordingly we denote the decomposition of $\Q$ by $\Q=\Q^+\sqcup\Q^0\sqcup\Q^-$.
For $s<0$ the typical fibres of $\Z^s\to M$ are hyperboloids
of two sheets, which are mutually identified via the central projection.
Similarly for $s>0$, where these are hyperboloids of one sheet.
Hence for any $s>0$ and $s<0$, the $s$-twistor space $\Z^s$ is identified
with the projectivization $\Proj\Q^+$ and $\Proj\Q^-$, respectively,
and we use the notation
$\Z^\pm:=\Z^{\pm1}\cong\Proj\Q^\pm$.
However, for $s=0$,  the typical fibre is the cone of null-vectors.
Hence 
$\Z^0=\Q^0$
and its projectivization
is a circle bundle over $M$ which will play a distinguished role later.
Altogether, we consider just three types of $s$-twistor spaces distinguished by the sign of $s$. 

Now, the almost para-quaternionic structure induces an almost $\eps$-complex structures on the respective twistor spaces.
The following statement is formulated as Proposition 6 in \cite{Alekseevsky2008a}.

\begin{prop}[\cite{Alekseevsky2008a}]    \label{p:stwistor}
  Let $(M,\Q)$ be an almost para\-/quaternionic manifold and let $\eps\in\{-1,0,1\}$.
  Any para\-/quaternionic connection induces a natural almost
  $\eps$-complex structure $\J^\eps$ on the $\eps$-twistor space $\Z^\eps$.
\end{prop}

The construction works roughly as follows.
A para\-/quaternionic connection $\nabla$ gives rise to a horizontal distribution $H^\nabla\subset T\Z^\eps$, complementary to the vertical subbundle of the projection $p:\Z^\eps\to M$.
The vertical subspace at any $z\in\Z^\eps$ is identified with the tangent space of an appropriate quadric in the oriented Minkowski space, hence it carries a canonical $\eps$-complex structure.
Next, any $z\in\Z^\eps$ is by definition an almost $\eps$-complex structure in $T_{p(z)}M$, and this lifts up to $H^\nabla_z$ via the inverse map of $T_zp$.
The two pieces then assembles into a natural almost $\eps$-complex structure on $T_z\Z^\eps$.
Moreover, for a section $s:M\to\Z^\eps$ of the projection $p$, let us denote by $J^s$ the corresponding almost $\eps$-complex structure on $M$.
Then the following compatibility relation holds:
\begin{equation} \label{eq:TpJTs}
J^s=Tp\o\J^\eps\o Ts.
\end{equation}

It is a natural question when two para-quaternionic connections induce the same $\eps$-complex structure on  the $\eps$-twistor space.
For $\eps=\pm 1$, this is carefully studied in \cite{David2009} and \cite{Ivanov2010}.
In particular, it turns out that all minimal para-quaternionic connections induce the same almost $\eps$-complex structure.
Such structure is therefore called \textit{canonical}.

The main outcome of the previous construction is that the integrability of the almost para\-/quaternionic structure is fully controlled by the integrability of the canonical almost $(\pm1)$-complex structures on twistor spaces.
The following statement is extracted from Theorem~21 in \cite{David2009}.

\begin{thm}[\cite{David2009}] \label{t:paraint}
  Let $(M,\Q)$ be an almost para\-/quaternionic manifold of dimension $2n>4$.
  Let $(\Z^\eps,\J^\eps)$ be the $\eps$-twistor space with the canonical almost $\eps$-complex structure, where $\eps=\pm1$.
  Then $\Q$ is integrable if and only if $\J^\eps$ is integrable.
\end{thm}
We revise this statement in section \ref{intgrblt}, where we also offer an extension to the case $\eps=0$.

Integrability of a $\eps$-complex structure  is equivalent to the vanishing of the corresponding Nijenhuis tensor.
Given a smooth manifold $Z$ and an endomorphism $A\in\End(TZ)=\Om^1(Z,TZ)$, the \textit{Nijenhuis tensor} of $A$ is given by the Fr\"olicher--Nijenhuis bracket, $N_A:=\frac12[A,A]\in\Om^2(Z,TZ)$, 
i.e.
\begin{equation} \label{eq:nijenhuis}
N_A(\xi,\eta):=-A^2[\xi,\eta]-[A\xi,A\eta]+A[A\xi,\eta]+A[\xi,A\eta],
\end{equation}
for any $\xi,\eta\in TZ$, where all brackets in \eqref{eq:nijenhuis} are the Lie brackets of vector fields.

  Note that for almost para-complex structures, 
  the Nijenhuis tensor vanishes if and only if the corresponding
  distributions are integrable in the sense of Frobenius.
  However, for almost tangent structures, vanishing of the Nijenhuis tensor is stronger than the integrability of the corresponding distribution,
  cf. e.g. \cite{Kobayashi1962}.

%%%
\section{Almost Grassmannian structures}  \label{Grassmann}
Here we collect several views on almost Grassmannian structures.
In particular, we emphasize the presence of the normal Cartan connection.
Recommended classical references are \cite{Bailey1991}, \cite{Akivis1996}, \cite{Machida2000} and \cite{Cap2009}.

\subsection{Grassmannians}\label{grassmn}
The Grassmannian of $p$-dimensional subspaces in $(p+q)$-dimensional real vector space, denoted as $\Gr_p(\R^{p+q})$, forms the model Grassmannian structure of type $(p,q)$.
The tangent space at each $\la\in\Gr_p(\R^{p+q})$ is naturally identified with the space of linear maps from $\la$ to the factor space $\R^{p+q}/\la$, i.e. with the tensor product $\la^*\otimes(\R^{p+q}/\la)$ of vector spaces of  dimensions $p$ and $q$.

Throughout this paper we consider just the structures of type $(2,n)$, where $n\ge2$.

The Grassmannian $\Gr_2(\R^{2+n})$ is the homogeneous space with the obvious transitive action of the Lie group $G:=PGL(2+n,\R)$, the quotient of the general linear group by its center 
(which consists of all real multiples of the identity). 
In particular, it coincides with $SL(2+n,\R)$, for odd $n$, and with $SL(2+n,\R)/\{\pm1\}$, for even $n$.
Denoting by $P$ the stabilizer of a 2-dimensional subspace in $\R^{2+n}$, we have $\Gr_2(\R^{2+n})\cong G/P$.
If this subspace is $\span{e_1,e_2}$, the span of the first two vectors of the standard basis of $\R^{2+n}$, then $P$ is represented by the block triangular matrices 
\begin{equation*}
  \PPmat{A\ &\ Z\\0\ &\ B}
\end{equation*}
with the blocks of sizes 2 and $n$ along the diagonal.
The subgroup $P\subset G$ is parabolic.

The related grading of the Lie algebra $\g=\sl(2+n,\R)$ is displayed in the following block form
\begin{equation*}
  \PPmat{\g_0\ &\ \g_1\\\g_{-1}\ &\ \g_0},
\end{equation*}
in particular, $\g_{-1}\cong\R^{2*}\otimes\R^n$,
$\g_0\cong\frak{s}(\gl(2,\R)\oplus\gl(n,\R))$ and $\g_1\cong\R^2\otimes\R^{n*}$.
The Lie algebra of $P$ is the sum of the nilpotent ideal $\g_1$ and the reductive subalgebra $\g_0$, $\p=\g_0\oplus\g_1$.
The central part of $\g_0$ consists of all multiples of the grading element, the semisimple part $\g_0^{ss}$ is isomorphic to the direct sum $\sl(2,\R)\oplus\sl(n,\R)$.

Let $G_0$ be the Lie subgroup in $P$ with the Lie algebra $\g_0$.
Evidently, 
\begin{equation}\label{eq:G0}
G_0\cong GL(2,\R)\. GL(n,\R),
\end{equation}
the quotient of the direct product of general linear groups by the subgroup consisting of all real multiples of the identity.
The adjoint representation of $G$ on $\g$ restricts to an injective group homomorphism $\Ad:G_0\to GL(\g_{-1})$.
Concretely, the action is given by
\begin{equation}\label{eq:Ad}
\Ad_{(A,B)}(X)=B\o X\o A^{-1},
\end{equation}
where the pair $(A,B)\in GL(2,\R)\x GL(n,\R)$ is a representative of an element of $G_0$ and $X\in\g_{-1}$ is seen as a linear map $\R^2\to\R^n$.

\subsection{Almost Grassmannian structures}  \label{grass}
An \textit{almost Grassmannian structure of type $(2,n)$} on a smooth manifold $M$ of dimension $2n\ge 4$ is given by an identification of the tangent bundle with the tensor product of two auxiliary vector bundles, 
\begin{equation}\label{eq:EF}
E^*\otimes F\xrightarrow{\cong}TM,
\end{equation}
where $\rank E=2$ and $\rank F=n$.

Note that an additional identification $\La^2E^*\cong\La^nF$ 
(or, equivalently, a trivialisation of the line bundle $\La^2 E\otimes\La^n F$)
is often taken as a part of the definition.
This just brings the notion of orientation into play;
the corresponding geometric structure is called  \textit{oriented} almost Grassmannian structure.
The model in this case is the  Grassmannian of oriented $2$-dimensional subspaces in $\R^{2+n}$.

Almost Grassmannian structures are G-structures with structure group as in \eqref{eq:G0}.
The structure group for the oriented version is the lift $S(GL(2,\R)\x GL(n,\R))\subset GL(2,\R)\x GL(n,\R)$ consisting of the indicated block matrices with determinant one.

A compatible connection of an almost Grassmannian structure is a linear connection on $TM\cong E\otimes F$ which is the tensor product of two linear connections on the auxiliary vector bundles $E$ and $F$.
An almost Grassmannian structure is called \textit{Grassmannian}, or \textit{integrable}, if there is a compatible torsion-free connection.
A natural class of normalized compatible connections is described in section \ref{weyl}.

Almost Grassmannian structures may be studied via the associated Segre structure, i.e. a field of Segre cones, see \cite{Akivis1996}, \cite{Grossman2000}, \cite{Mettler2013}.
Under the isomorphism \eqref{eq:EF}, the \textit{Segre cone} in $T_xM$, $x\in M$, is exactly the set of simple elements of $E_x^*\otimes F_x$, i.e. the set of linear maps $E_x\to F_x$ of rank one.
The Segre cone is doubly ruled by linear subspaces of dimensions 2 and $n$; 
the corresponding subbundles in $\Gr_2(TM)$ and $\Gr_n(TM)$ are denoted as  $\Cal F$ and $\Cal E$ and their elements are called  $\al$- and \textit{$\be$-planes}, respectively.
The notation reflects the fact that these subbundles are naturally identified with the projectivized auxiliary bundles so that $\Cal F\cong\Proj F$ and $\Cal E\cong\Proj E$.
An almost Grassmannian structure is called \textit{$\be$-integrable}, if any $\be$-plane from $\Cal E$ is tangent to a unique immersed $n$-dimensional submanifold of $M$ whose all tangent spaces are elements of $\Cal E$.
The notion of \textit{$\al$-integrability} is analogous.

\subsection{Normal Cartan connection} \label{grasscartan}
Throughout this paper we rely on the fact that almost Grassmannian structures can be described as parabolic Cartan geometries.
In particular, we have a canonical normalization condition determining a distinguished class of compatible connections.
We have to recall some generalities first.

In this paragraph, $G$ may denote an arbitrary Lie group, $P\subset G$ its Lie subgroup and $\p\subset\g$ the corresponding Lie algebras.
The model Cartan geometry associated to the homogeneous space $G/P$ consists of the homogeneous principal $P$-bundle $G\to G/P$ and the Maurer--Cartan form  $\om\in\Om^1(G,\g)$.
General Cartan geometry of type $G/P$ on a smooth manifold $M$ consists of a principal $P$-bundle $\G\to M$ and a {Cartan connection} $\om\in\Om^1(\G,\g)$.
In particular, $\om$ is an absolute parallelism, i.e. it provides a global identification $TG\cong G\x\g$.
Among other identifications determined by $\om$, the most frequent one is 
\begin{equation*}
TM\cong\G\x_P(\g/\p), 
\end{equation*}
where the right hand side reads as the associate bundle to $\G$ with the typical fibre $\g/\p$ and the natural action of $P$ (i.e. the one induced by the adjoint action on $\g$).
The \textit{curvature} of the Cartan geometry $(\G\to M,\om)$ is an element of $\Om^2(\G,\g)$ defined by 
\begin{equation*}
\ka:=d\om+\om\wedge\om.
\end{equation*}
Since the curvature is strictly horizontal, the corresponding frame form reduces to a $P$-equivariant map
$\G\to\La^2(\g/\p)^*\otimes\g$, the so-called curvature function.
Composing  with the quotient projection $\g\to\g/\p$, we obtain a $P$-equivariant map $\G\to\La^2(\g/\p)^*\otimes(\g/\p)$ representing a tensor field $\tau\in\Om^2(M,TM)$, which is called the \textit{torsion} of the Cartan geometry.

In the case that $G$ is semisimple and $P$ parabolic, the corresponding Cartan geometry is called \textit{parabolic}.
The pair $P\subset G$ from section \ref{grassmn} related to Grassmannians is of this type and, moreover, 
the length of the corresponding grading $\g=\g_{-1}\oplus\g_0\oplus\g_1$ is the smallest possible.
Parabolic geometries with this property are called \textit{$|1|$-graded}.
In contrast to general parabolic geometries, a lot of things simplifies if the structure is $|1|$-graded.
In the following we repeatedly enjoy this fact.

It turns out that (as for most $|1|$-graded parabolic geometries) the nilpotent subalgebra $\g_1\subset\p$ coincides with the first prolongation of $\g_0\subset\p$ and its second prolongation vanishes.
For $\g_{-1}\cong\R^{m}$, $m=\dim M$, the map from \eqref{eq:del} can be seen as the $G_0$-equivariant map
$\g_{-1}^*\otimes\g_0\to\La^2\g_{-1}^*\otimes\g_{-1}$ deduced (according to the gradation of $\g$) from the differential $\del$ in the chain complex computing the cohomology of the Lie algebra $\g_{-1}$ with coefficients in $\g$.
In this context, we also use the duality between $\g/\p\cong\g_{-1}$ and $\g_1$ via the Cartan--Killing form.
In particular, the curvature function may be seen as a map $\G\to\La^2\g_1\otimes\g$.

The natural normalization condition is given by the $P$-equivariant map $\del^*$, the codifferential  in the complex computing the Lie algebra homology of $\g_1$ with coefficients in $\g$:
the parabolic geometry is called \textit{normal} if its curvature function takes values in $\ker\del^*\subset\La^2\g_1\otimes\g$.
In such case, the composition with the quotient projection $\ker\del^*\to\ker\del^*/\im\del^*$ yields a new quantity, the \textit{harmonic curvature}.
It follows that harmonic curvature determines the full curvature and has an interpretation in underlying terms.

In fact, the maps $\del^*$ and $\del$ are adjoint with respect to an appropriate inner product.
This gives rise to the $G_0$-equivariant self-adjoint endomorphism $\square:=\del\o\del^*+\del^*\o\del$, the so-called Kostant Laplacian, which determines a Hodge decomposition of the chain complex.
In particular, the kernel of this operator,
$$
\ker\square\subset\ker\del^*\subset\La^2\g_1\otimes\g,
$$
is isomorphic to the second homology group.
It follows that the lowest non-zero homogeneous component of the curvature function has values in $\ker\square$, i.e. it coincides with the corresponding homogeneous component of the harmonic curvature.
The nice thing is that  $\ker\square$ is algorithmically computable as a $G_0$-representation.

The following statement is the starting point of our further considerations, cf. \cite[sec.~2--3]{Machida2000} and \cite[sec.~4.1.3]{Cap2009}:

\begin{prop}[\cite{Machida2000}, \cite{Cap2009}] \label{p:grasscartan}
An almost Grassmannian structure of type $(2,n)$ on $M$ is equivalent to a normal parabolic geometry $(\G\to M,\om)$ of type $G/P$, where $G=PGL(2+n,\R)$ and $P$ is the parabolic subgroup as above.
In terms of \eqref{eq:EF}, the components of the harmonic curvature are indicated in the following tables:

\medskip
\begin{tabular}{cc}
\parbox[c][][c]{0.45\textwidth}{
\begin{center}
$n=2$ \\[6pt] 
\begin{tabular}{|c|l|}
  \hline homog. & {section of} \\
  \hline 2 & {$S^2E\otimes\La^2F^*\otimes\sl(E)$}\\
  \hline 2 & {$\La^2E\otimes S^2F^*\otimes\sl(F)$}\\
  \hline
\end{tabular}
\end{center}}&
\parbox[c][][c]{0.45\textwidth}{
\begin{center}
$n>2$ \\[6pt]
\begin{tabular}{|c|l|}
  \hline homog. & {section of} \\
  \hline 1 & {$S^2E\otimes\La^2F^*\otimes E^*\otimes F$}\\
  \hline 2 & {$\La^2E\otimes S^2F^*\otimes\sl(F)$}\\
  \hline
\end{tabular}
\end{center}}
\end{tabular}
\medskip

In particular, the torsion of the Cartan geometry vanishes for $n=2$ and coincides with the harmonic curvature component of homogeneity one for $n>2$.
\end{prop}

More details on the indicated decomposition and the torsion component are in sections \ref{decomp} and \ref{connections}, respectively.

Note that the notion of $\be$-, respectively $\al$-integrability of the almost Grassmannian structure is controlled by the vanishing of the first, respectively second component in the displayed tables, see the results of 
\cite{Goncharov1987}, \cite{Akivis1999} and \cite{Machida2000}.

In the case $n=2$, it follows that the almost Grassmannian structure is equivalent to a conformal structure of split signature so the two harmonic curvatures correspond to anti-self-dual, respectively self-dual, part of the Weyl curvature tensor.%
\footnote{Be aware that the conventions in references are not always consistent;
we follow the one in which the $\be$-integrability corresponds to the anti-self-duality.}
More comments on this special case are in section \ref{dim4}.

\subsection{Weyl connections} \label{weyl} 
As for any parabolic geometry, there is a natural class of compatible connections, the \textit{Weyl connections}.
These are in a bijective correspondence with reductions of the principal $P$-bundle $\G\to M$ to the structure group $G_0\subset P$.
Equivalently, they correspond to global $G_0$-equivariant sections of the canonical projection $\G\to\G_0$, where $\G_0\to M$ is the quotient principal bundle with structure group $G_0$, i.e. $\G_0=\G/\exp\g_1$ in the current setting.
The family of Weyl connections is parametrized by one-forms on $M$.
The difference between two such connections, $\nabla$ and $\wh\nabla$, is expressed via $\Y\in\Om^1(M)$ so that
\begin{equation}
\wh\nabla_{\xi}{\eta} =\nabla_{\xi}{\eta}-\{\{\Y,\xi\},\eta\},
\label{eq:hatna}
\end{equation}
for any $\xi,\eta\in\Ga(TM)$.
Each bracket on the right hand side is the algebraic bracket induced by the one in the Lie algebra $\g=\g_{-1}\oplus\g_0\oplus\g_1$, see \cite[sec.~5.1.6]{Cap2009}.
In particular, $\{\Y,\xi\}$ is an endomorphism of $TM\cong E^*\otimes F$, pointwise corresponding to elements of $\g_0$.

As for any $|1|$-graded parabolic geometry, 
all  Weyl connections share the same torsion, namely, the torsion of the Cartan connection $\om$.
In these terms, the integrability of an almost Grassmannian structure is equivalent to the vanishing of the torsion of the corresponding normal Cartan connection, i.e. to the vanishing of the harmonic curvature component of homogeneity one.
This condition is automatically satisfied for $n=2$.

\subsection{Grassmannian circles}\label{circles}
Belonging to the broad family of parabolic geometries, almost Grassmannian structures admit classes of distinguished curves in the sense of \cite{Cap2004}.%
\footnote{Besides an enormous terminology related to concrete geometries, 
these curves have various general nicknames, e.g.  Cartan's circles, generalized geodesics or canonical curves.}
According to the absolute parallelism $T\G\cong\G\x\g$ determined by $\om$, they are given as projections of flow lines of constant vector fields corresponding to elements from $\g_{-1}$.
Equivalently, they are the curves that develop to the orbits of one-parameter groups in the homogeneous model with generators in $\g_{-1}$.
We note that many geometric properties of such curves are controlled by the algebraic properties of the map $\g_1\to\g_{-1}$ given by $Z\mapsto [[Z,X],X]$, where $X\in\g_{-1}$ is a representing element of a respective type of curve and the brackets are the Lie brackets in $\g$.

It follows that, for almost Grassmannian structures of type $(2,n)$, there are two types of distinguished curves:
those that are everywhere, respectively nowhere, tangent to the Segre cone; in this article we call them \textit{null}, respectively \textit{generic, Grassmannian circles}.
The former curves are given by a tangent vector in one point so that a collinear tangent vector yields just different parametrization of the same path.
In particular, null Grassmannian circles are common unparametrized geodesic of all compatible connections.
The latter curves are given by an initial condition of second order.
It also follows that any Grassmannian circle is a geodesic of a Weyl connection satisfying some additional condition.
This point of view is applied in section \ref{curves}.

\subsection{Twistor correspondence for almost Grassmannian structures}
\label{grasstwist}
Here we describe the twistor spaces of \cite{Bailey1991}, respectively \cite{Machida2000}, following the vocabulary of \cite{Cap2005} and \cite{Cap2009}.
Let $G=PGL(2+n,\R)$ and $P\subset G$ be as above.
Let $P'\subset G$ be the stabilizer of the line $\span{e_1}$ 
spanned by the first vector of the standard basis in $\R^{2+n}$.
Hence $Q:=P\cap P'$ is the stabilizer of the flag
$\span{e_1}\subset\span{e_1,e_2}$.
Alike $G/P$ was identified with the Grassmannian $\Gr_2(\R^{2+n})$, the homogeneous space  $G/P'$ is identified
with the projective space $\R\P^{1+n}$ and $G/Q$ with the proper flag manifold.
The flag manifold $G/Q$ is fibred both over the Grassmannian $G/P$ and over
the projective space $G/P'$.
It is called the \textit{correspondence space} of $G/P$ and $G/P'$, 
while the latter spaces are its \textit{twistor spaces}.
For later use, we figure the respective subgroups of $G$ in the block matrix
form:
\begin{equation*}
  \begin{gathered}
    Q=\Pmat{a&b&Z_1\\0&d&Z_2\\0&0&B},\\ 
    P=\Pmat{a&b&Z_1\\c&d&Z_2\\0&0&B},\quad
    P'=\Pmat{a&b&Z_1\\0&d&Z_2\\0&Y&B},
  \end{gathered}
\end{equation*}
where the separators distinguish the blocks of sizes 2 and $n$ as before; in particular, $a,b,c,d\in\R$.
Note that all these subgroups are parabolic.

Let $(\G\to M,\om)$ be the normal parabolic geometry of type $G/P$ associated to an almost Grassmannian structure on $M$.
The \textit{correspondence space} of $M$ with respect to $Q\subset P$ is the orbit space
$$
\Cal CM:=\G/Q,
$$
the total space of the fibre bundle over $M$ whose typical fibre is $P/Q\cong\R\P^1$.
It easily follows that elements of $\Cal CM$ correspond to 1-dimensional subspaces in the rank 2 auxiliary vector bundle from \eqref{eq:EF}, i.e. $\Cal CM\cong\Proj E$.

The restricted Cartan geometry $(\G\to\Cal CM,\om)$ is a parabolic geometry of type $G/Q$, which is automatically normal, but not necessarily regular.
The regularity means that all homogeneous components of the curvature function have positive degree.
(This condition is satisfied trivially for $|1|$-graded geometries; the current length of gradation corresponding to the Lie subalgebra $\q\subset\g$ of the parabolic subgroup $Q\subset G$ is two.)

Regular and normal parabolic geometries of type $G/Q$ are equivalent to the so-called \textit{generalized path geometries}.
Such structure on $\Cal CM$ consists of two subbundles $D,V\subset T\Cal CM$ of rank 1 and $n$, respectively, with trivial intersection and some other properties; see \cite{Grossman2000} or \cite[sec.~4.4.3]{Cap2009}.
Under the identification $T\Cal CM\cong\G\x_Q(\g/\q)$, the two subbundles in $T\Cal CM$ are 
\begin{equation} \label{eq:DV}
  D\cong\G\x_Q(\p/\q),\quad V\cong\G\x_Q(\p'/\q),
\end{equation}
where $\q$, $\p$ and $\p'$ is the Lie algebra to $Q$, $P$ and $P'$,
respectively.
Clearly, the line subbundle $D\subset T\Cal CM$ is the vertical subbundle of the projection $\Cal CM\to M$.

The integrability of the almost Grassmannian structure is reflected on the correspondence space level as follows, see \cite[Prop.~4.4.5]{Cap2009}.

\begin{prop}[\cite{Cap2009}]  \label{p:4.4.5}
Let $E^*\otimes F\xrightarrow{\cong}TM$ be an almost Grassmannian structure of type $(2,n)$ on $M$, let $(\G\to M,\om)$ be the corresponding normal parabolic geometry and let $(\G\to\Cal CM,\om)$ be the normal parabolic geometry over the correspondence space.
Then the former parabolic geometry is torsion-free (i.e. the almost Grassmannian structure on $M$ is integrable) if and only if the letter parabolic geometry is regular.
\end{prop}

Note that the stated property is automatically satisfied in the case $n=2$.

Generalized path geometry generalizes the notion of \textit{path geometry}, which is a system of unparametrized curves that are determined by a tangent direction in one point.
A path geometry on a smooth manifold $X$ induces a generalized path geometry on the projectivized tangent bundle $\Proj TX$ so that the paths on $X$ coincides with the projections of the integral curves of the distribution $D$.
The complementary distribution $V$ corresponds to the vertical subbundle of the projection $\Proj TX\to X$, in particular, it is involutive.
We return to this topic in section \ref{ident}.

%%%%
\section{Equivalence and first interactions} \label{Inter}
Here we come with first couple of interactions.
Most of them are expected from the previous preparations and they should not be surprising; 
we only make these expectations precise.
The main observations of this section may be summarized as follows:

\begin{thm}\label{p:first}
Let $M$ be a smooth manifold of dimension $2n\ge4$.
There is a natural bijective correspondence between (equivalent classes of) 
almost para\-/quaternionic structures $\Q\subset\End(TM)$ and
almost Grassmannian structures $E^*\otimes F\xrightarrow{\cong}TM$ of type $(2,n)$.
Under this identification:
\begin{enumerate}[(1)]
\item
The eigenspaces of para-complex structures (equivalently, the kernels of almost tangent structures) from $\Q\subset\End(TM)$ are just the $\be$-planes, the maximal linear subspaces contained in the Segre cone of $E^*\otimes F\cong TM$.
\item
Minimal para-quaternionic connections are just the Weyl connections of the associated normal Cartan connection.
\item 
A null $Q$-planar curve is a null Grassmannian circle if and only if it is a common unparametrized geodesic of all compatible connections.
\item 
A generic $Q$-planar curve is a Grassmannian circle if and only if the equation \eqref{eq:S} is satisfied.
\end{enumerate}
\end{thm}

\subsection{Equivalence of structures} \label{paragrass}
Both almost para\-/quaternionic structures and almost Grassmannian structures of type $(2,n)$ can be regarded as G-structures.
We have vaguely referred to the structure group of the former structure in section \ref{minimal}, while the structure group of the latter structure is described in section \ref{grassmn}.
Passing to the vector space level $T_xM\cong\R^{2n}$, $x\in M$, this is a subgroup of $GL(T_xM)\cong GL(2n,\R)$ up to a covering.
At this stage, the equivalence of the two structures is easy to see, cf. e.g.  \cite[sec.~4.3]{Alekseevsky2008a}:

Consider we are given a vector space $W$ of dimension $2n$ and a linear isomorphism $W\cong\R^{2*}\otimes\R^n$.
Then an endomorphism of $\R^2$ gives rise to an endomorphism of $W$ via the action on the first factor.
The restriction just to the trace-free endomorphisms of $\R^2$ yields a 3-dimensional subspace of endomorphisms of $W$.
This obviously defines a para\-/quaternionic structure, which we call the \textit{standard para\-/quaternionic structure} and denote by $\Q_{std}$.
Conversely, given a para\-/quaternionic structure $\Q_W$ on $W$, the algebra $\span{\id}+\Q_W$ is isomorphic to the algebra of para\-/quaternions, i.e. to the matrix algebra $Mat_{2\x2}(\R)$.
Any irreducible $Mat_{2\x2}(\R)$-module is isomorphic to $\R^2$, hence the  $Mat_{2\x2}(\R)$-module $W$ is isomorphic to the tensor product $\R^2\otimes\R^n$.
Under this identification, the action of $\Q_W$ on $W$ corresponds to the action on the first factor.

More concretely, let $X$ be an element of $\R^{2*}\otimes\R^n$, seen as a linear map $X:\R^2\to\R^n$, and let $\un{A}$ be a trace-free endomorphism of $\R^2$.
Then the corresponding element $A$ of $\Q_{std}$, i.e. an endomorphism of $\R^{2*}\otimes\R^n$, is given by
\begin{equation} \label{eq:compos}
A(X)=X\o\underline{A}.
\end{equation}
In these terms, the norm squared on $\Q_{std}$, defined by \eqref{eq:A^2}, corresponds to the determinant, 
\begin{equation}
|A|^2=\det\underline{A}.
\label{eq:det}
\end{equation}

Now, the interpretation of the structure group of an almost Grassmannian structure in terms of the corresponding para-quaternionic structure may be seen as follows.
According to section \ref{grassmn}, we have $\R^{2*}\otimes\R^n\cong\g_{-1}$ and the restricted adjoint representation identifies the structure group $G_0$ with a subgroup of $GL(\g_{-1})\cong GL(2n,\R)$.
Let an element of $G_0$ be represented by a pair $(C,D)\in GL(2,\R)\x GL(n,\R)$, let $f=\Ad_{(C,D)}\in GL(\g_{-1})$ be the corresponding linear isomorphism and let an element $A\in\Q_{std}$ be represented by $\un{A}\in\End(\R^2)$.
Then \eqref{eq:Ad} and \eqref{eq:compos} yield
\begin{equation*}
f(A(X))
=D\o(X\o\un{A})\o C^{-1}
=(D\o X\o C^{-1})\o(C\o\un{A}\o C^{-1})
=\phi(A)(f(X)),
\end{equation*}
where $\phi(A)$ denotes the element of $\Q_{std}$ corresponding to $C\o\un{A}\o C^{-1}\in\End(\R^2)$.
Thus, $f=\Ad_{(C,D)}$ is a para-quaternionic automorphism of the standard para-quaternionic structure on $\g_{-1}$.
Conversely, if $f\in GL(\g_{-1})$ is such an automorphism, then the defining condition \eqref{eq:auto} translates under the current notation as
\begin{equation*}
f(A(X))=f(X)\o\un{\phi}(\un{A}),
\end{equation*}
where $\un{\phi}$ is an algebra automorphism of $\End(\R^2)$, i.e. an automorphism of the algebra of para-quaternions.
From section \ref{paraquat} we know that $\un{\phi}$ can be represented by an element $C\in SL(2,\R)$ so that $\un{\phi}(\un{A})=C\o \un{A}\o C^{-1}$. 
It is now easy to see that $f=\Ad_{(C,D)}$, where $D\in GL(n,\R)$ is determined by the condition $D\o X=f(X)\o C$.
To summarize,

\begin{lem}
A para\-/quaternionic structure $\Q_W$ on a vector space $W$ of dimension $2n$ is equivalent to an isomorphism $W\cong\R^{2*}\otimes\R^n$ so that $\Q_W$ corresponds to the standard para\-/quaternionic structure.
Under this identification, the Lie group $G_0$ from \eqref{eq:G0} coincides with the group of para-quaternionic automorphisms of $\Q_W$.
\end{lem}

Passing to the Lie algebra level, 
$\g_0$ is the direct sum of the semisimple part $\sl(2,\R)\oplus\sl(n,\R)$ and one-dimensional center.
This allows the interpretation $\g_0\cong\sl(2,\R)\oplus\gl(n,\R)$.
The first summand consists of trace-free endomorphisms of $\R^2$, i.e. of elements of the standard para-quaternionic structure on $\g_{-1}\cong\R^{2*}\otimes\R^n$.
The second summand consists of all endomorphisms of $\R^n$, i.e. of endomorphisms of $\g_{-1}$ commuting with $\Q_{\g_{-1}}$.
In particular, $\gl(n,\R)$ and $\g_0$ may be seen as the centralizer and the normalizer, respectively, of $\Q_{\g_{-1}}\cong\sl(2,\R)$ in $\gl(\g_{-1})\cong\gl(2n,\R)$.
This point of view is employed in \cite{David2009}, see also section \ref{connections}.
 
\medskip

Now we can characterize the Segre cone in $\R^{2*}\otimes\R^n$ as follows.
By definition,  a linear map $X:\R^2\to\R^n$ belongs to the Segre cone if and only if the kernel of $X$ has dimension one.
Arbitrary complementary subspace $\ell$ to $\ker X$ in $\R^2$ determines a para-complex structure $\underline{A}$ so that $\ell$ and $\ker X$ is its eigenspace corresponding to 1 and $-1$, respectively.
Then,  according to \eqref{eq:compos}, $X$ is an eigenvector (corresponding to the eigenvalue 1) of the associated para-complex structure $A$ on $\R^{2*}\otimes\R^n$.
Conversely, let $X$ be an $(+1)$-eigenvector of a para-complex structure $A\in\Q_{std}$ and let $\underline{A}$ be the corresponding para-complex structure on $\R^2$; i.e. $A(X)=X\o\underline{A}=X$.
Since $\underline{A}$ is not the identity, $X$ cannot be of full rank and, since $X\ne0$, it has rank one.
Hence $X$ belongs to the Segre cone.
Alternatively, an eigenspace of a para-complex structure can be realized as the kernel (image) of some tangent structure, and vice versa.
Altogether,

\begin{lem}
With previous identifications,
a non-zero element of $W\cong\R^{2*}\otimes\R^n$ belongs to the Segre cone if and only if it is an eigenvector of a para-complex structure (equivalently, lies in the kernel of a tangent structure) from $\Q_W\cong\Q_{std}$.
\end{lem}

As a consequence, we note that the eigenspaces of para-complex structures (equivalently, the kernels of tangent structures) form the maximal linear subspaces contained in the Segre cone.
Altogether, the first part of Theorem \ref{p:first} follows.

\subsection{Decompositions}  \label{decomp}
The decomposition of complex forms into $(p,q)$-types has the following counterparts.
Suppose we are given a real vector space $W$ endowed with an endomorphism $A\in\End(W)$ which squares to a multiple of the identity, written $A^2=A\o A=-|A|^2\id$ as in \eqref{eq:A^2}.
Let us consider the bilinear maps $\ph:W\x W\to W$.
The notion of the \textit{type $(p,q)$ of $\ph$ with respect to $A$} is given as follows:
\begin{itemize}
  \item $\ph$ is of type $(1,1)$ if $\ph(AX,AY)=-A^2\ph(X,Y)=|A|^2\ph(X,Y)$,
  \item $\ph$ is of type $(0,2)$ if $\ph(AX,Y)=\ph(X,AY)=-A\ph(X,Y)$,
  \item $\ph$ is of type $(2,0)$ if  $\ph(AX,Y)=\ph(X,AY)=A\ph(X,Y)$,
\end{itemize}
where, here and after, all identities are meant to hold for all $X,Y\in W$.

If $|A|^2\ne 0$ then $\ph$ decomposes uniquely into the sum of components of
particular types with respect to $A$, namely,
$\ph=\ph^{2,0}_A+\ph^{1,1}_A+\ph^{0,2}_A$, 
where
\begin{equation}  \label{eq:pqpart}
  \begin{aligned}
    \ph^{1,1}_A(X,Y)&:=\frac1{2|A|^2}\left(|A|^2\ph(X,Y)+\ph(AX,AY)\right),\\
    \ph^{0,2}_A(X,Y)&:=\frac1{4|A|^2}\left(|A|^2\ph(X,Y)-\ph(AX,AY)+A\ph(AX,Y)+A\ph(X,AY)\right),\\
    \ph^{2,0}_A(X,Y)&:=\frac1{4|A|^2}\left(|A|^2\ph(X,Y)-\ph(AX,AY)-A\ph(AX,Y)-A\ph(X,AY)\right).
  \end{aligned}
\end{equation}
The respective elements $\ph^{p,q}_A$ are called the \textit{$(p,q)$-parts of $\ph$ with respect to $A$}.

If $|A|^2=0$, the notion of type of $\ph$ is rather degenerate.
E.g., there are forms which are simultaneously of all three types with respect to $A$ so any type decomposition is a priori meaningless.
However, in order to unify the later treatment and simplify some formulations, we will use the label \textit{$(0,2)$-part of $\ph$ with respect to $A$} also for the map given by
\begin{equation}  \label{eq:02part}
\ph^{0,2}_A(X,Y):=\frac14\left(-\ph(AX,AY)+A\ph(AX,Y)+A\ph(X,AY)\right),
\quad \text{where}\ |A|^2=0.
\end{equation}

If $W$ carries a para\-/quaternionic structure $\Q_W\subset\End(W)$, we may consider bilinear maps $W\x W\to W$ that are of previous types with respect to all of $A\in\Q_W$.
Moreover, the notion of type $(1,1)$ has a good meaning also  for bilinear forms $W\x W\to\R$.
The space of bilinear forms of type $(1,1)$ with respect to all $A\in\Q_W$ is denoted by $\bigotimes^{1,1}W^*$.
There is a natural projection 
$\bigotimes^2W^*\to\bigotimes^{1,1}W^*$, $\ph\mapsto\ph^{1,1}$,
given by 
\begin{equation} \label{eq:pi11} 
  \ph^{1,1}(X,Y):=\frac14\left(\ph(X,Y)-\ph(IX,IY)-\ph(JX,JY)+\ph(KX,KY)\right),
\end{equation}
where $(I,J,K)$ is an arbitrary basis of $\Q_W$, i.e. a triple satisfying
\eqref{eq:IJK}.
The kernel of this projection is a natural complementary subspace to
$\bigotimes^{1,1} W^*$ in $\bigotimes^2 W^*$.
In particular, this yields the decomposition 
\begin{equation}  \label{eq:kerpi}
  \La^2 W^*=\La^{1,1} W^*\oplus\ker\pi^{1,1},
\end{equation}
where $\pi^{1,1}$ denotes the restriction of the projection above to
$\La^2 W^*\subset\bigotimes^2 W^*$.

\smallskip

From the previous subsection we know that any para\-/quaternionic structure on a $2n$-dimensional vector space $W$ can always be viewed as the standard para-quaternionic structure under an identification $W\cong \R^{2*}\otimes \R^n$. 
This yields the  decomposition
\begin{equation} \label{eq:la2}
  \La^2 W^*\cong (\La^2\R^2\otimes S^2\R^{n*}) \oplus (S^2\R^2\otimes\La^2\R^{n*}).
\end{equation}
In order that our endeavour have some meaning, this must agree with the decomposition in \eqref{eq:kerpi}.

\begin{lem} \label{l:decomp}
  With the current notation, the following hold:
  $$
  \La^2\R^2\otimes S^2\R^{n*}=\La^{1,1}(\R^2\otimes \R^{n*})
  \qtext{and}
  S^2\R^2\otimes\La^2\R^{n*}=\ker\pi^{1,1}.
  $$
\end{lem}
\begin{proof}
The standard para-quaternionic structure on $\R^{2*}\otimes \R^{n}$ is given by the action on $\R^2$.
For any $A\in\Q_{std}$ and any $\ph\in\La^2\R^2\otimes S^2\R^{n*}$, we have $\ph(AX,AY)=\det\un{A}\.\ph(X,Y)$, where $\un{A}$ is the corresponding endomorphism of $\R^2$.
According to \eqref{eq:det}, we see that $\ph$ is of type $(1,1)$ so $\La^2\R^2\otimes S^2\R^{n*}\subseteq\La^{1,1}(\R^2\otimes \R^{n*})$.

Any element of $S^2\R^2\otimes\La^2\R^{n*}$ is a linear combination of simple elements $e_i\odot e_j\otimes v^k\wedge v^l$, where $(e_i)$ is the standard basis of $\R^2$ and $(v^i)$ is the standard basis of $\R^{n*}$.
Let us choose a basis $(I,J,K)$ of $\Q_{std}$ so that it corresponds to the matrices as in \eqref{eq:ijk}.
According to \eqref{eq:pi11}, it follows that 
\begin{equation*}
\pi^{1,1}(e_1\odot e_2\otimes v^k\wedge v^l)
=\frac14(e_1\odot e_2+e_1\odot e_2-e_2\odot e_1-e_2\odot e_1)\otimes v^k\wedge v^l
=0
\end{equation*}
and, similarly, that $\pi^{1,1}(e_1\odot e_1\otimes v^k\wedge v^l)=\pi^{1,1}(e_2\odot e_2\otimes v^k\wedge v^l)=0$, for any $v^k,v^l\in\R^{n*}$.
Thus, $S^2\R^2\otimes\La^2\R^{n*}\subseteq\ker\pi^{1,1}$.

Now, the statement follows from the  complementarity of the respective subspaces, i.e. from \eqref{eq:kerpi} and \eqref{eq:la2}.
\end{proof}

\subsection{Distinguished connections}\label{connections}
Both minimal para-quaternionic connections and Weyl connections of a normal Cartan connection are affine connections that are compatible with the structure in question and share the same normalized torsion.
The respective normalization conditions, i.e. the invariant subbundles in $\La^2T^*M\otimes TM$, are described in two different ways according to either para-quaternionic or Grassmannian (parabolic) terminology.
In previous two subsections we explained the equivalence of the two geometric structures, and this was used for a double expression of the decomposition of the space of 2-forms.
With a bit finer discussion we can directly show the two normalization conditions  coincide.
We recall that there is no torsion in dimension four so the only non-trivial discussion concerns the general case corresponding to $n>2$.

On the one hand, the normalization condition from \cite[sec.~3]{David2009} is described in terms of the $(p,q)$-type decompositions as follows.
Firstly, for a basis $(I,J,K)$ of $\Q_W$, let the endomorphism of $\La^2 W^*\otimes W$ be defined by 
\begin{equation*}
\Pi(\ph):=\frac23\left(\ph^{0,2}_I+\ph^{0,2}_J+\ph^{0,2}_K\right),
\end{equation*}
where the individual summands are as in \eqref{eq:pqpart}.
It turns out that the definition is independent of the basis of $\Q_W$, $\Pi$ is a projector (i.e. $\Pi\o \Pi=\Pi$),
and its kernel coincides with the image of $\del:W^*\otimes\gl(n,\R)\to\La^2 W^*\otimes W$.
Here $\gl(n,\R)$ denotes the centralizer of $\Q_{std}\cong\sl(2,\R)$ in $\gl(W)\cong\gl(2n,\R)$ as discussed in section \ref{paragrass}.
Thus, the image of $\Pi$ is a complementary subspace to the image of $\del$.
In particular, it is contained in the kernel of $\pi^{1,1}$.

Secondly, extending the Lie algebra to $\sl(2,\R)\oplus\gl(n,\R)\cong\g_0$, the natural $G_0$-invariant complement of the image of $\del:W^*\otimes\g_0\to\La^2 W^*\otimes W$ is described within the image of $\Pi$ as
\begin{equation*}
\Cal C:=\left\{\ph\in\im \Pi : \tr(A\o\ph(X,-))=0,\ \text{for all\ $X\in W$ and $A\in\{I,J,K\}$}\right\}.
\end{equation*}

On the other hand, for $W=\g_{-1}\subset\g$, the natural complement of the image of $\del:\g_{-1}^*\otimes\g_0\to\La^2\g_{-1}^*\otimes\g_{-1}$ is given by the kernel of $\del^*$, namely, $\Cal D:=(\ker\del^*)\cap(\La^2\g_{-1}^*\otimes\g_{-1})$.
Since $\g$ is $|1|$-graded,
this subspace is actually harmonic, i.e. it coincides with the  corresponding irreducible component in the kernel of  the Kostant Laplacian, cf. section \ref{grasscartan}.
Since $\g_{-1}\cong\R^{2*}\otimes\R^n$, the space $\La^2\g_{-1}^*\otimes\g_{-1}$ decomposes according to \eqref{eq:la2}.
In the summand corresponding to $\ker\pi^{1,1}$ we have two obvious traces, 
\begin{equation*}
\R^2\otimes\La^2\R^{n*}\otimes\R^n 
\xleftarrow{\tr_1}
S^2\R^2\otimes\La^2\R^{n*}\otimes\R^{2*}\otimes\R^n 
\xrightarrow{\tr_2}
S^2\R^2\otimes\R^{n*}\otimes\R^{2*}.
\end{equation*}
It is shown in \cite[sec.~4.1.3]{Cap2009} that $\Cal D$ is characterized as the intersection of the kernels of these two traces,
\begin{equation}
\Cal D=\left\{\ph\in S^2\R^2\otimes\La^2\R^{n*}\otimes\R^{2*}\otimes\R^n : \tr_1(\ph)=\tr_2(\ph)=0 \right\}.
\label{eq:D}
\end{equation}

\begin{lem}\label{l:CD}
The two normalization conditions, i.e. the invariant subspaces $\Cal C,\Cal D\subset\La^2W^*\otimes W$, coincide.
\end{lem}

\begin{proof}
As a typical element of $\Cal D$ we choose
\begin{equation}
\ph:=e_1\odot e_1\otimes \ups^1\wedge\ups^2\otimes\ep^2\otimes u_3,
\label{eq:phi}
\end{equation}
where $(e_i)$ denotes the standard basis of $\R^2$, $(\ep^i)$ its dual basis of $\R^{2*}$,
$(u_i)$ denotes the standard basis of $\R^n$ and $(\ups^i)$ its dual basis of $\R^{n*}$.
We also choose a basis $(I,J,K)$ of $\Q_{std}$ so that the corresponding  matrices are as in \eqref{eq:ijk}.
Now, it is an easy exercise to show that 
$\ph^{0,2}_I=\ph$ and $\ph^{0,2}_J+\ph^{0,2}_K=\frac12\ph$.
Hence $\Pi(\ph)=\ph$, i.e. $\ph\in\im\Pi$.
As a consequence of $\tr_2(\ph)=0$, we also see that the full trace of $A\o\ph(X,-)$ vanishes for any $X\in\R^{2*}\otimes\R^n$ and $A\in\Q_{std}$.
Thus, $\ph\in\Cal C$ and, since $\Cal D$ is an irreducible $G_0$-representation, it follows that $\Cal D\subseteq\Cal C$.
Since $\Cal D$ and $\Cal C$ are both complementary to the same subspace, we have $\Cal D=\Cal C$.
\end{proof}

Altogether, the second part of Theorem \ref{p:first} follows. 

Note that, as representations of the reductive Lie group $G_0$, the source and the target space of the map $\del$ are completely reducible. 
A closer look on the decompositions shows that all components on both sides appear with multiplicity one.
From this, and the fact that $\del$ is either trivial or an isomorphism on each irreducible component, it follows that the invariant complement to $\im\del$ is actually unique.
From this perspective the previous lemma is no surprise.

\smallskip
In conclusion, we describe the difference between two compatible connections having the same torsion.
On the one hand, two Weyl connections differ as in \eqref{eq:hatna}.
Expanding the difference term in our case yields
\begin{equation}
\{\{\Y,\xi\},\eta\}
=-\Y(\xi)\eta-\Y(\eta)\xi
=-\xi\o\Y\o\eta-\eta\o\Y\o\xi,
\label{eq:Ups}
\end{equation}
where the second expression is due to the interpretation of $\xi,\eta$, respectively $\Y$, as fields of linear maps $E\to F$, respectively $F\to E$.
On the other hand, the difference of two para-quaternionic connections is displayed in \eqref{eq:ups}.
This and \eqref{eq:Ups} have to be just two distinct expressions of the same difference term.
For a concrete choice of basis $(I,J,K)$, following the identifications from section \ref{paragrass}, one easily shows they indeed coincide (up to a constant multiple).

\subsection{Distinguished curves}\label{curves}
$Q$-planar curves of an almost para-quaternionic structure are defined in terms of compatible affine connections, see section \ref{Q-planar}.
Grassmannian circles of an almost Grassmannian structure are defined via the associated Cartan connection, see section \ref{circles}.
We start with several easy observations on $Q$-planar curves, then we remind an alternative characterization of Grassmannian circles in terms of compatible affine connections, which will lead to a comparison of these two families of distinguished curves.

Let $\Q\subset\End(TM)$ be an almost para-quaternionic structure and let $TM\cong E^*\otimes F$ be the corresponding almost Grassmannian structure. 
Let $\ga$ be a curve and let us interpret the tangent vector field $\dot\ga$ as a field of linear maps $E\to F$ along $\ga$.
According to the reasoning in section \ref{paragrass}, the $Q$-planarity of $\ga$ with respect to a compatible connection $\nabla$ may be written as
\begin{equation}
\nabla_{\dot\ga}\dot\ga=\dot\ga\o{S},
\label{eq:Qplan}
\end{equation}
where ${S}\in\End(E)$.
Expressing the $Q$-planarity of $\ga$ with respect to another compatible connection, it follows from \eqref{eq:Ups} that the corresponding endomorphism  changes as
\begin{equation}
\wh{S}={S}+2\Y\o\dot\ga,
\label{eq:hatS}
\end{equation}
where $\Y$ is interpreted as a field of linear maps $F\to E$ as above.
As we already observed before, geodesics of a compatible connection are $Q$-planar with respect to this (equivalently, to any) compatible connection.
For the converse we have to distinguish two cases:

If $\ga$ is generic, we see that $\Y$ may always be chosen so that \eqref{eq:hatS} vanishes, i.e. so that $\ga$ is a geodesic of the corresponding connection.
If $\ga$ is null, this is not the case. 
But if $\ga$ is a geodesic of one compatible connection, then it is a (unparametrized) geodesic of all of them
(this is because $\dot\ga\o\Y\o\dot\ga$ is a multiple of $\dot\ga$, for any $\Y$).
Note that we may indeed refer to any compatible connection, since the torsion does not play any role concerning geodesics.
Altogether, we summarize as

\begin{lem}
Let $\ga$ be a $Q$-planar curve of some (equivalently, any) compatible connection, seen as an unparametrized curve.
\begin{enumerate}[(1)]
\item
If $\ga$ is generic, then it is a  geodesic of some compatible connection.
\item
If $\ga$ is null, then it is a  geodesic of some compatible connection if and only if it is a geodesic of any compatible connection.
\end{enumerate}
\end{lem}

Concerning Grassmannian circles,
their alternative definition is provided by the so-called Rho (or Schouten) tensor, which is associated to any Weyl connection of any parabolic geometry.
The following formulations are specialized to the $|1|$-graded case.
Let $\si:\G_0\to\G$ be the $G_0$-equivariant section corresponding to a Weyl connection $\nabla$ and let $\om_1\in\Om^1(\G,\g_1)$ be the $\g_1$-part of the Cartan connection $\om$.
The \textit{Rho tensor} $\Rho$ is given by the pullback $\si^*\om_1\in\Om^1(\G_0,\g_1)$.
By the horizontality of $\om$ and the identification $\g_1\cong\g_{-1}^*$, it may be seen as a tensor field $\Rho\in\Om^1(M,T^*M)$.
The Rho tensor transforms under the change of Weyl connection as
\begin{equation}
\wh\Rho(\xi)
=\Rho(\xi)+\nabla_\xi\Y+\frac12\{\Y,\{\Y,\xi\}\}
=\Rho(\xi)+\nabla_\xi\Y-\Y\o\xi\o\Y,
\label{eq:hatRho}
\end{equation}
for any $\xi\in\Ga(TM)$, where we follow that same conventions as in \eqref{eq:Ups}, see \cite[sec.~5.1.8]{Cap2009}.

Now, as a special case of a more general setting \cite[sec.~5.3.1]{Cap2009}, a curve $\ga$ is a Grassmannian circle if and only if there is a Weyl connection $\nabla$ such that $\ga$ is its geodesic and the corresponding Rho tensor $\Rho$ vanishes for the tangent vectors of $\ga$, i.e.
\begin{equation}\label{eq:narho}
\nabla_{\dot\ga}\dot\ga=0
\qtext{and}
\Rho(\dot\ga)=0.
\end{equation}
From this characterization, it is clear that Grassmannian circles are $Q$-planar curves of $\nabla$ and hence of all compatible connections.
The identification of former class of curves among the latter one is as follows:

If $\ga$ is null, then from general properties of null Grassmannian circles and the previous lemma we see that
a null $Q$-planar curve is a null Grassmannian circle if and only it is a common unparametrized geodesic of all compatible connections.
Hence the third part of Theorem \ref{p:first} follows.

If $\ga$ is generic, then the description of Grassmannian circle among the $Q$-planar curves can be easily adapted from \cite[sec.~6.2]{Bailey1991} to our setting.
Namely, Theorem 6.4 of that reference translates to a characterization in terms of the invariant differential equation,%
\footnote{Comparing with the original formulation we differ in the sign in the front of the term containing $\Rho$.
This just reflects the difference in the definition of Rho tensor here, which we took from \cite{Cap2009}, and in \cite{Bailey1991}.}
\begin{equation}\label{eq:S}
\nabla_{\dot\ga}S=\frac12 S\o S+2\Rho(\dot\ga)\o\dot\ga,
\end{equation}
where $\nabla$ is a compatible connection, $S$ is an endomorphism given by \eqref{eq:Qplan}, $\Rho$ is the Rho tensor of $\nabla$ and $\Rho(\dot\ga)$ is seen as a field of linear maps $F\to E$.
The invariance of this equation follows from the transformation formulas for $\nabla$, $S$ and $\Rho$ with respect to a change of compatible connection, cf. \eqref{eq:Ups}, \eqref{eq:hatS} and \eqref{eq:hatRho}.
Clearly, a curve satisfying \eqref{eq:narho} satisfies \eqref{eq:S}, i.e. a Grassmannian circle is a $Q$-planar curve satisfying the latter equation.
The converse statement is shown by suitable changes of compatible connections.
Altogether, the fourth part of Theorem \ref{p:first} follows.

%%%%%
\section{Twistor spaces revised}\label{Twistor}
In this section we recover the canonical almost $\eps$-complex structures on $\eps$-twistor spaces and the corresponding integrability statement from section \ref{paratwist}, which we extend also to the case $\eps=0$.
Then  we comment the 0-twistor space in detail, both in general and integrable case.
This provides a link between the two notions of twistor correspondence for almost para\-/quaternionic and almost Grassmannian structures.

\subsection{Setup}  \label{setup}
Let $(M,\Q)$ be an almost para\-/quaternionic manifold, equivalently an almost Grassmannian structure, 
and let $(\G\to M,\om)$ be the induced normal parabolic geometry of type
$G/P$.
For any $x\in M$, the Cartan connection $\om$ identifies $(T_xM,\Q_x)$ with $(\g_{-1},\Q_{std})$, where $\Q_{std}$ denotes the standard para\-/quaternionic structure on $\g_{-1}\cong\R^{2*}\otimes\R^n$.
This is represented by trace-free endomorphisms of $\R^2$, which form the  Lie algebra $\sl(2,\R)$, and this is seen as the left-upper block from the matrix description of $\g^{ss}_0$ in section \ref{grassmn}.

The parabolic subgroup $P$ acts on $\Q_{std}$ via the adjoint action so that the orbits of the action consist of those elements which have the same norm.
For any $\eps\in\{-1,0,1\}$, let us choose an $\eps$-complex structure $j^\eps\in\Q_{std}$ and let us denote by $R^\eps\subset P$ the stabilizer of $j^\eps$.
In other words, $R^\eps$ is the subgroup consisting of all para\-/quaternionic automorphisms of $(\g_{-1},\Q_{std})$ which commute with $j^\eps$.
Hence each orbit is the homogeneous space $P/R^\eps$ and this is the typical fibre of the $\eps$-twistor bundle $\Z^\eps\to M$ defined in section \ref{paratwist}.
Hence 
\begin{equation}  \label{eq:Z}
\Z^\eps\cong\G\x_P(P/R^\eps)\cong\G/R^\eps.
\end{equation}
and the Cartan geometry $(\G\to M,\om)$ gives rise to a Cartan geometry $(\G\to\Z^\eps,\om)$ of type $G/R^\eps$ on each $\eps$-twistor space.
Note that none of these Cartan geometries is parabolic.
Nevertheless, the Cartan connection $\om$ provides the identification
\begin{equation*}
T\Z^\eps\cong\G\x_{R^\eps}(\g/\r^\eps), 
\label{eq:TZ}
\end{equation*}
where $\r^\eps$ is the Lie algebra to $R^\eps$.
To summarize,

\begin{lem}
For each $\eps\in\{-1,0,1\}$, the $\eps$-twistor space is $\Z^\eps\cong\G/R^\eps$ and it carries a canonical Cartan geometry $(\G\to\Z^\eps,\om)$ of type $G/R^\eps$.
\end{lem}

Obviously, only the semisimple part of the left-upper block of $P$ acts non-trivially on $\Q_{std}$.
Thus, the typical fibres of $\eps$-twistor bundles may be identified as
\begin{equation*}
P/R^{-}\cong SL(2,\R)/SO(2),\quad
P/R^{0}\cong SL(2,\R)/\R_+,\quad
P/R^{+}\cong SL(2,\R)/SO(1,1),
\end{equation*}
where $\R_+$ stands for the additive group of real numbers, which is realized as the subgroup of the form $\pmat{1&b\\0&1}$ in $SL(2,\R)$.

Considering the $G_0$-principal bundle $\G_0=\G/\exp\g_1$ as in section \ref{weyl} and $R^\eps_0:=G_0\cap R^\eps$, 
the identification \eqref{eq:Z} can be written as $\Z^\eps\cong\G_0/R^\eps_0$, cf. \cite{Alekseevsky2008a}, \cite{Mettler2013}.
The individual subgroups $R^\eps_0$ are isomorphic to
\begin{equation*}
R^-_0\cong SO(2)\.GL(n,\R),\quad
R^0_0\cong \R_+\.GL(n,\R),\quad
R^+_0\cong SO(1,1)\.GL(n,\R).
\end{equation*}

In concrete computations we use the $\eps$-complex structures $j^\eps\in\Q_{std}$, whose $2\x 2$-blocks in the previously indicated matrix description are as follows:
\begin{equation} \label{eq:j}
j^-=\pmat{0&1\\-1&0},\quad 
j^0=\pmat{0&1\\0&0},\quad 
j^+=\pmat{1&0\\0&-1}.
\end{equation}
(Be aware of a small abuse of notation, which also applies below.)
The explicit description of the corresponding subgroups $R^\eps\subset P$ and their Lie algebras $\r^\eps\subset\p$ yields that elements of $\g/\r^\eps$ may be represented by the matrices of the form
\begin{equation*}
\Bmat{b&c&0\\c&-b&0\\X_1&X_2&0}\in\g/\r^-,\quad
\Bmat{b&0&0\\e&-b&0\\X_1&X_2&0}\in\g/\r^0,\quad
\Bmat{b&d&0\\-d&-b&0\\X_1&X_2&0}\in\g/\r^+.
\end{equation*}

%%%%
\subsection{Induced $\eps$-complex structures} \label{sicomplex}
The current point of view allows an alternative description of almost $\eps$-complex structures on $\eps$-twistor spaces. 
In contrast to the development in section \ref{paratwist}, they are now described via the associated Cartan connection and related identifications.

\begin{prop}    
Let $p:\Z^\eps\to M$ be the $\eps$-twistor space, $\eps\in\{-1,0,1\}$, of an almost para\-/quaternionic manifold $(M,\Q)$.
Then the total space $\Z^\eps$ carries a unique almost $\eps$-complex structure $\J^\eps$ such that \eqref{eq:TpJTs} holds.
\end{prop}
\begin{proof}
For each $\eps\in\{-1,0,1\}$, let $j^\eps$, $R^\eps$ and $\r^\eps$ be as above.
Let us define an endomorphism $J^\eps:\g/\r^\eps\to\g/\r^\eps$ by 
  \begin{equation} \label{eq:J}
    \BBmat{U&*\\X&*}\mapsto \BBmat{Uj^\eps&*\\Xj^\eps&*}
         =\BBmat{U&*\\X&*}\.\PPmat{j^\eps&0\\0&0} \mod\r^\eps.
  \end{equation}
  Obviously, $J^\eps\o J^\eps=\eps\id$ and it is easy to check that
  $J^\eps$ is also $R^\eps$-invariant.
  Hence it gives rise to an almost $\eps$-complex structure $\J^\eps$ on
  $\Z^\eps$.

By the identification \eqref{eq:Z}, a section $s$ of the projection $\Z^\eps\to M$ is represented by a $P$-equivariant function $\si:\G\to P$ such that $s(x)=u\si(u) R^{\eps}$, for each $x\in M$ and any $u\in\G_x$.
The tangent vector $\xi\in T_xM$ is represented by the couple $[u,X+\p] \in\G\x_P(\g/\p)$.
Remember that any other representative of the same equivalence class is of the form $[u p,\Ad_{p^{-1}}X+\p]$, for some $p\in P$.
In these terms, the action of the $\eps$-complex structure $J^s$ on $T_xM$ corresponding to $s$  is given by
  \begin{equation}  \label{eq:J^s}
    [u,X+\p]\mapsto[u,X\.\Ad_{\si(u)}j^\eps+\p],
  \end{equation}
Next, the tangent map to the section $s:M\to\Z^\eps$ is written as
\begin{equation*}
[u,X+\p]\mapsto[u\si(u),\Ad_{\si(u)^{-1}}X+\r^\eps], 
\end{equation*}
whereas the tangent map to the projection $p:\Z^\eps\to M$ is just $[u,X+\r^\eps]\mapsto[u,X+\p]$.
Altogether, the composition $Tp\o\J^\eps\o Ts$ maps
\begin{equation*}
[u,X+\p]\mapsto[u\si(u),\Ad_{\si(u)^{-1}}X\.j^\eps+\p],
\end{equation*}
which  clearly coincides with \eqref{eq:J^s}.
Thus, $Tp\o\J^\eps\o Ts=J^s$ and the equality \eqref{eq:TpJTs} holds.

A direct computation shows that if $J$ is an $R^\eps$-invariant $\eps$-complex structure on $\g/\r^\eps$ then $J$ coincides with $J^\eps$ up to the sign  in the case $\eps=\pm1$, respectively up to a non-zero real multiple in the case $\eps=0$.
However, the previous paragraph reveals that the condition \eqref{eq:TpJTs} is satisfied if and only if $J$ and $J^\eps$ coincide. 
  Hence the almost $\eps$-complex structure $\J^\eps$ on $\Z^\eps$ is unique.
\end{proof}

Of course, this almost $\eps$-complex structure must recover the canonical one from section \ref{paratwist}, which was determined by a (arbitrary) minimal para-quaternionic connection.
The relation can be made explicit with the following observations.
From section \ref{connections} we know that the minimal para\-/quaternionic connections are exactly the Weyl connections of the associated normal Cartan connection.
Any such connection $\nabla$ is given by a reduction $\G_0\to M$ of the Cartan $P$-bundle $\G\to M$ to $G_0\subset P$.
Hence $TM\cong\G_0\x_{G_0}(\g/\p)$ and $T\Z^\eps\cong\G_0\x_{R^\eps_0}(\g/\r^\eps)$, where $R^\eps_0=G_0\cap R^\eps$ as before.
In  this description, the vertical subbundle of the projection $p:\Z^\eps\to M$ corresponds to the subspace $\p/\r^\eps\subset\g/\r^\eps$, while the horizontal subbundle $H^\nabla\subset T\Z^\eps$ corresponds to $(\g_-\oplus\r^\eps)/\r^\eps$, the unique subspace in $\g/\r^\eps$ that is both $R^\eps_0$-invariant and complementary to $\p/\r^\eps$.
Now, the original description of $\J^\eps$ in terms of its horizontal and vertical part can be readily compared with the current invariant approach.

\subsection{Integrability} \label{intgrblt}
In the setting of section \ref{setup}, 
the torsion of the Cartan geometry $(\G\to M,\om)$ and $(\G\to\Z^\eps,\om)$
is denoted by $\tau\in\Om^2(M,TM)$ and $\T\in\Om^2(\Z^\eps,T\Z^\eps)$, respectively.
By definitions, $\T$ is strictly horizontal with respect to the projection
$p:\Z^\eps\to M$, hence 
$\tau(\xi,\eta)=Tp(\T(\hat \xi,\hat \eta))$,
where $\hat \xi,\hat \eta\in T\Z^\eps$ are any lifts of $\xi,\eta\in TM$.
In other words,
\begin{equation}
  \tau=Tp\o\T\o(Ts\x Ts)
  \label{eq:TpTauTsTs}
\end{equation}
for any section $s:M\to\Z^\eps$.

The following lemma can be seen as a Cartan-geometric analogue of the well-known fact
that an almost (para-)complex structure is integrable if and only if the
$(0,2)$-part of the torsion of some (and consequently any) compatible affine
connection vanishes. 
The reasoning below is very similar to the one  in \cite[sec.~4.4.10]{Cap2009}.
An alternative treatment in the case $\eps=1$ can be found in \cite[sec.~5]{Alekseevsky2008a}.

\begin{lem} \label{l:intgrblt}	
  Let $\Z^\eps$ be the $\eps$-twistor space with the canonical almost $\eps$-complex structure, $\eps\in\{-1,0,1\}$, and let $\T$ be the torsion of the associated Cartan connection over $\Z^\eps$.
  Then the Nijenhuis tensor of $\J^\eps$ is a non-zero constant multiple of the $(0,2)$-part of $\T$ with respect to $\J^\eps$, which is taken according to the definition in \eqref{eq:pqpart}, respectively \eqref{eq:02part}.
\end{lem}
\begin{proof}
  To deal efficiently with the tensor fields on $\Z^\eps$ we use the corresponding frame forms
  with respect to the Cartan connection $\om$. 
  On the one hand,
  the frame form of the torsion $\T$ is the 
  $R^\eps$-equivariant functions
  $\G\to\La^2(\g/\r^\eps)^*\otimes(\g/\r^\eps)$,
  which assigns to each $u\in\G$ the bilinear map
  \begin{equation*}
    (X+\r^\eps,Y+\r^\eps)\mapsto
    \pi\left([X,Y]-\om([\omi(X)(u),\omi(Y)(u)])\right),
  \end{equation*}
  where $\pi$ is the quotient projection $\g\to\g/\r^\eps$.
  Similarly, the frame form of $\J^\eps$ is the constant function
  $\G\to(\g/\r^\eps)^*\otimes(\g/\r^\eps)$ with value $J^\eps$, which is
  described in \eqref{eq:J}.
  Now one can express the frame form of $\T^{0,2}_{\J^\eps}$ following 
  the conventions from section \ref{decomp}, distinguishing the cases $\eps\ne0$
  and $\eps=0$.

  On the other hand, the frame form of the Nijenhuis tensor $N_{\J^\eps}$, cf. \eqref{eq:nijenhuis}, is the equivariant function, which assigns to each $u\in\G$ the bilinear
  map
  \begin{equation*}
    \begin{aligned}
      (X+\r^\eps&,Y+\r^\eps) \mapsto  \\
      -&(J^\eps)^2(\pi(\om([\omi(X)(u),\omi(Y)(u)])))
      -\pi(\om([\omi(J^\eps X)(u),\omi(J^\eps Y)(u)]))+\\
      +&J^\eps(\pi(\om([\omi(J^\eps X)(u),\omi(Y)(u)])))
      +J^\eps(\pi(\om([\omi(X)(u),\omi(J^\eps Y)(u)]))),
    \end{aligned}
  \end{equation*}
  where $J^\eps X$ denotes any element in $\g$ such that $\pi(J^\eps X)=J^\eps(\pi(X))$.
  Note that by the $R^\eps$-invariancy of the $\eps$-complex structure $J^\eps$, this is indeed a well-defined object.
  
  Let us consider the tensor field  
  $\Cal S:=N_{\J^\eps}-4\T^{0,2}_{\J^\eps}$.
  Taking into accounts that $(J^\eps)^2=\eps\id$, a simple substitution shows the frame form of $\Cal S$ is the constant function assigning to each $u\in\G$ the bilinear map
  \begin{equation*}
    \begin{aligned}
      (X+\r^\eps&,Y+\r^\eps) \mapsto \\
      -&(J^\eps)^2(\pi([X,Y]))-\pi([J^\eps X,J^\eps Y])
      +J^\eps(\pi([J^\eps X,Y]))+J^\eps(\pi([X,J^\eps Y])).
    \end{aligned}
  \end{equation*}
  However, by the definition of $J^\eps$ in \eqref{eq:J}, it immediately follows
  that 
  $$
  (J^\eps)^2(\pi([X,Y]))=J^\eps(\pi([X,J^\eps Y]))
  \quad\text{and}\quad
  \pi([J^\eps X,J^\eps Y])=J^\eps(\pi([J^\eps X,Y])),
  $$
  for any $\eps\in\{-1,0,1\}$.
  Therefore $\Cal S=0$, which completes the proof.
\end{proof}

Here is the promised extension and reinterpretation of the statement cited as Theorem \ref{t:paraint}.
An analogous statement in the four-dimensional case is formulated in section \ref{dim4}.

\begin{thm} \label{p:intgrblt}	
  Let $(M,\Q)$ be an almost para\-/quaternionic manifold of dimension $2n>4$.
  Let $(\Z^\eps,\J^\eps)$ be the $\eps$-twistor space with the canonical almost $\eps$-complex structure, where $\eps\in\{-1,0,1\}$.
  Then $\Q$ is integrable if and only if $\J^\eps$ is integrable.
\end{thm}

\begin{proof}
On the one hand, the integrability of $\Q$ is equivalent to the vanishing of the torsion $\tau$ of the associated normal Cartan connection $\om$ over $M$, which equals to the harmonic curvature component of homogeneity one, cf. Proposition \ref{p:grasscartan}.
On the other hand, the integrability of $\J^\eps$ is equivalent to the vanishing of the $(0,2)$-part of the torsion $\T$ of $\om$ understood as a Cartan connection over $\Z^\eps$, cf. Lemma \ref{l:intgrblt}.

Let $\Q$ be integrable, i.e. $\tau=0$.
Hence the whole Cartan curvature is determined by the harmonic curvature component of homogeneity two. 
By the description in Proposition \ref{p:grasscartan}, its frame form takes values in $\sl(n,\R)$, the lower right block in $\g_0$, and the curvature component of homogeneity three has necessarily values in $\g_1$.
Hence, for each $\eps$, the Cartan curvature takes values  in $\r^\eps$.
This means that the torsion $\T$ is also trivial, i.e. $\J^\eps$ is integrable. 

Conversely, let $\J^\eps$ be integrable, i.e.  the torsion $\T$ has trivial $(0,2)$-part with respect to $\J^\eps$. 
The main consequence of the relations \eqref{eq:TpJTs} and \eqref{eq:TpTauTsTs} is that $\tau$ must have vanishing $(0,2)$-part with respect to any $\eps$-complex structure contained in $\Q$.
Since $\tau$ coincides with the harmonic torsion, it is heavily restricted.
Namely, its frame form takes values in the irreducible $G_0$-representation $\Cal D$ as described in \eqref{eq:D}.
For each $\eps$, we are going to show that there is actually no non-zero element of $\Cal D$ that would satisfy the requirement.
Hence this requirement can be satisfied if and only if $\tau=0$, which is equivalent to the integrability of $\Q$.

Since $\Cal D$ is an irreducible representation, it is enough to find, for each $\eps$, a concrete $\eps$-complex structure with respect to which an arbitrarily chosen element of $\Cal D$ has non-vanishing $(0,2)$-part.
As a representative element $\ph\in\Cal D$ we choose the one as in \eqref{eq:phi}
and, for $\eps\in\{-1,0,1\}$, we choose the elements $j^\eps\in\Q_{std}$ as in \eqref{eq:j}.
Now it easily follows that $\ph$ is of type $(0,2)$ with respect to $j^+$, 
has non-trivial $(0,2)$-part with respect to  $j^-$ and also with respect to $j^0$.
\end{proof}

\subsection{The 0-twistor space} \label{ident}
We have number of fibre bundles over $M$, arising either from the para-quaternionic or the Grassmannian side.
On the one hand, the bundle $\Q\subset\End(TM)$ defining an almost para\-/quaternionic structure on $M$ is decomposed into disjoint subbundles such that $\Q^0\subset\Q$ coincides with the 0-twistor space $\Z^0$, see section \ref{paratwist}, 
On the other hand, for the corresponding almost Grassmannian structure $E^*\otimes F\cong TM$, the projectivization $\Proj E$ is identified with the bundle $\Cal E$ of $\be$-planes in $TM$
and this is further identified with the correspondence space $\Cal CM$, see sections \ref{grass} and \ref{grasstwist}.
The following statement supplies yet another identification that relates the two universes.

\begin{prop}   \label{p:z0}
  Let $\Z^0$ be the 0-twistor space over an almost para\-/quaternionic manifold $M$, let $\Cal CM$ be the correspondence space of the corresponding almost Grassmannian structure on $M$ and let $D\subset T\Cal CM$ be the vertical subbundle of the projection $\Cal CM\to M$.
Then $\Z^0$ is naturally identified with $D$.
In particular, $\Proj\Z^0$ is identified with $\Cal CM$. 
\end{prop}
\begin{proof} 
  As an associated bundle over $\Cal CM$, the line bundle $D$ is identified
  with $\G\x_Q(\p/\q)$, see \eqref{eq:DV}. 
  The action of $Q$ on $\p/\q$ is transitive and the stabilizer
  of any element is just the subgroup $R^0$ from section \ref{setup}.
  Hence $D\cong\G\x_Q(Q/R^0)\cong\G/R^0\cong\Z^0$, according to
  \eqref{eq:Z}.
\end{proof}

As a consequence of Proposition \ref{p:4.4.5}, an integrable Grassmannian structure on $M$ gives rise to a generalized path geometry on $\Cal{C}M$.
The one-dimensional distribution $D\subset T\Cal CM$ is just the vertical subbundle of the projection $\Cal CM\to M$.
If $n>2$, 
the distribution $V\subset T\Cal CM$ is automatically involutive, which allows to construct a leaf space $X$
so that $\Cal{C}M$ is locally identified with $\Proj TX$ and $V$ corresponds to the vertical subbundle of the projection $\Proj TX\to X$, see \cite[Prop.~4.4.4]{Cap2009}.
In particular, the generalized path geometry on $\Cal{C}M$ is locally equivalent to a path geometry on $X$ so that the points in $M$ corresponds to the paths in $X$.
In conclusion, we have an additional interpretation of the 0-twistor space in the integrable case:

\begin{prop} \label{p:z0int}
In addition to assumptions of the previous proposition, let $\dim M>4$, let the structure on $M$ be integrable and let $X$ be a local leaf space of the foliation determined by the involutive distribution $V\subset T\Cal CM$.
Then  $\Z^0$ is locally identified with the tangent bundle $TX$ so that the rank $n+1$ distribution $\ker\J^0\subset T\Z^0$,
with $\J^0$ being the canonical 0-complex structure on $\Z^0$,
corresponds to the vertical subbundle of the canonical projection $TX\to X$.
\end{prop}

\begin{proof}
Following \cite[Prop.~4.4.4]{Cap2009}, we recall some details on the local identification of $\Cal CM$ with $\Proj TX$.
Denoting by $\psi:\Cal CM\supset U\to X$ the local leaf space projection, its tangent map $T\psi$ induces a linear isomorphism $T_xU/V_x\to T_{\psi(x)}X$, for any $x\in U$. 
Hence $D_x\subset T_xU$ projects to a one-dimensional subspace in $T_{\psi(x)}X$, i.e. an element in $\Proj T_{\psi(x)}X$ which is denoted as $\tilde\psi(x)$. 
It is shown the tangent map to $\tilde\psi:U\to\Proj TX$ is invertible, therefore $\tilde\psi$ is an open embedding.
It is now easy to see that $\tilde\psi$ extends to a local embedding of $D$ into $TX$.

From Proposition \ref{ident} we know that $\Z^0$ coincides with $D$,
hence $\Z^0$ is locally identified with $TX$.
The projection $\Z^0\cong TX\to X$ factorizes through $\Proj TX\cong\Cal CM$ and we already know that the vertical subbundle of $\Proj TX\to X$ coincides with $V$.
It is enough to show  that, under the canonical projection $\Z^0\to\Proj\Z^0$, $\ker\J^0$ maps to $V$.

By the proof of Proposition \ref{sicomplex}, the almost complex structure
  $\J^0$ corresponds to the $R^0$-invariant endomorphism
  $J^0:\g/\r^0\to\g/\r^0$ given by \eqref{eq:J}.
  With the same conventions  as before, 
  the kernel of $J^0$ is the $R^0$-invariant subspace represented by the
  matrices of the form
  $$\Bmat{u&0&0\\0&-u&0\\0&X_2&0}.
  $$
  The tangent map to the canonical projection
  $\Z^0\cong\G/R^0\to\G/Q\cong\Proj\Z^0$ corresponds to the obvious
  $R^0$-invariant projection $\g/\r^0\to\g/\q$ determined by $\r^0\subset\q$.
  The image of $\ker J^0$ in $\g/\q$ is then represented by 
  $$\Bmat{0&0&0\\0&0&0\\0&X_2&0}.
  $$
  Now we see that the image coincides with the $Q$-invariant subspace
  $\p'/\q\subset\g/\q$, which defines the distribution $V\subset
  T\Cal CM$ as in \eqref{eq:DV}.
\end{proof}

%%%
\section{Remarks}  \label{Final}
Here we add two things: 
several necessary remarks on the four-dimensional case and a note on compatible metrics.

\subsection{Dimension four}  \label{dim4}
As we repeatedly noticed,  the case when the base manifold $M$ has dimension four (i.e. the case $n=2$ according to the previous notation) is quite specific.
While the four-dimensional para-quaternionic structures are sometimes considered as a degenerate case, it is well known that almost Grassmannian structures of type $(2,2)$ are equivalent to conformal structures of split signature.
In terms of distinguished directions in the tangent bundle $TM$, the relation is such that the Segre cone of $TM\cong E^*\otimes F$ is just the cone of the non-zero null-vectors of the conformal structure.
Note that the Segre cone forms a hyper-quadric in the tangent space exactly in this dimension.

On the level of Lie algebras, with the description as in section \ref{grassmn}, $\g=\frak{sl}(4,\R)$ and the block corresponding to $\g_{-1}$ is of size $2\x2$. 
Let us consider the quadratic form on $\g_{-1}$ defined by the determinant;
the corresponding polar form is denoted by $\de$ for later purposes.
Evidently, the null-vectors of this form exhaust exactly the Segre cone
of rank-one elements in $\g_{-1}\cong\R^{2*}\otimes\R^2$.
The adjoint action of $G_0$ on $\g_{-1}$ changes the form conformally, which
leads to the identification $G_0\cong CSO_0(2,2)$.
For oriented almost Grassmannian structures, the structure group is a
two-fold covering of the just mentioned one, namely, $G_0\cong CSpin(2,2)$.
Under this identification, the bundles $E$ and $F$ are identified with the two spinor bundles.
Hence the correspondence space $\Cal CM$, as defined in section \ref{grasstwist}, is identified with the projectivized spinor bundle.
The two harmonic curvature components from Proposition \ref{p:grasscartan} corresponds to the self-dual and the anti-self-dual part of the conformal Weyl curvature, cf. e.g. \cite[sec.~4.1.4]{Cap2009}.

Concerning another notions from section \ref{Grassmann}, we just remark that the Weyl connections for conformal structures are the torsion-free connections preserving the conformal class of metrics and that null, respectively generic, Grassmannian circles coincide with null geodesics, respectively conformal circles.

The development of section \ref{Inter} includes also the case $n=2$, only the discussion on the normalization condition in section \ref{connections} is vacuous as there is no torsion in that case.
The description of the para-quaternionic structure in terms of the conformal one is as follows.
It is an easy observation that the inner product $\de$ is, up to a non-zero constant multiple, the unique non-degenerate bilinear form on $\g_{-1}$ which is of type $(1,1)$ with respect to the standard para-quaternionic structure $\Q_{std}$.
This means that, for any $A\in\Q_{std}$ and $X,Y\in\g_{-1}$, the following holds:
$$
\de(AX,AY)=|A|^2\de(X,Y).
$$
If $|A|^2\ne 0$ then this condition is equivalent to 
$$
\de(AX,Y)+\de(X,AY)=0,
$$
i.e. $A$ is skew with respect to $\de$.
For $|A|^2=0$, the latter condition is stronger.
Conversely, it turns out that if $A$ is an endomorphism of $\g_{-1}$, which is skew with respect to $\de$ and whose square $A^2$ is a multiple of the identity, then $A$ belongs to $\Q_{std}$.
Altogether, we have a characterization of the standard para-quaternionic structure in terms of $\de$, which is obviously independent of a  multiple of $\de$.
The geometric interpretation of these observations is the following:
an endomorphism $A$ of the tangent bundle of a para-quaternionic 4-manifold $(M,\Q)$ belongs to $\Q\subset\End(TM)$ if and only if $A\o A$ is a multiple of the identity and $A$ is skew with respect to any metric form the conformal class of the corresponding conformal structure.

Finally, let us consider the $\eps$-twistor spaces with the canonical almost $\eps$-complex structures from section \ref{Twistor}.
Everything works fine for this dimension up to the following adjustment of Theorem \ref{p:intgrblt}
(the structure on $M$ is automatically integrable so this is no relevant condition).
According to the description of the harmonic curvatures in Proposition \ref{p:grasscartan}, it follows that vanishing of the first component in the corresponding table is a sufficient condition for the integrability of the induced almost $\eps$-complex structure.
That this condition is also necessary follows by the very same scenario as in the proof of Theorem \ref{p:intgrblt}.
By remarks after Proposition \ref{p:grasscartan}, this condition corresponds to the anti-self-duality of the corresponding conformal structure (respectively, to the $\be$-integrability of the Grassmannian structure).

Altogether, we conclude with
\begin{prop}
Let $(M,\Q)$ be a 4-dimensional para-quaternionic manifold, let $[g]$ be the corresponding conformal structure on $M$ and let $\eps\in\{-1,0,1\}$.
\begin{enumerate}[(1)]
\item The $\eps$-twistor space $\Z^\eps$ is identified with the subbundle of $\End(TM)$ of those elements which square to $\eps\id$ and which are skew with respect to $[g]$.
\item The canonical almost $\eps$-complex structure on $\Z^\eps$ is integrable if and only if the conformal structure on $M$ is anti-self-dual.
\end{enumerate}
\end{prop}

For $\eps=-1$, the characterization of the respective twistor space may be shortened by saying that $\Z^{-1}$ consists of orthogonal almost complex structures in $TM$.
This should commemorate the classical formulations, cf. \cite{Atiyah1978} and \cite[Prop.~4.4.11]{Cap2009}.

Note that the interpretation of the 0-twistor space as in Proposition \ref{p:z0int} has to be adjusted accordingly in this dimension.
I.e., the assumption of integrability has to be substituted by the anti-self-duality of the conformal structure (respectively, by the $\be$-integrability of the Grassmannian structure).
This is what is needed to form a local leaf space $X$, the rest remains the same.

\subsection{Compatible metrics}  \label{Metric}
It is a very important situation if there exists a (pseudo-)Riemannian metric which is compatible with the given geometric structure.
This is thoroughly studied both from the para-quaternionic and the Grassmannian point of view.
Following \cite{Alekseevsky2008a}, \cite{Alekseevsky2009} and \cite{Bailey1991}, let us quickly summarize some of the classical issues here.

There is a natural decomposition of the bundle $S^2T^*M$ in the spirit of section \ref{decomp}.
On the one hand, the almost para-quaternionic structure $\Q\subset\End(TM)$ induces the decomposition 
\begin{equation*}
  S^2T^*M=S^{1,1}T^*M\oplus \ker\pi^{1,1},
\end{equation*}
where $\pi^{1,1}:S^2T^*M\to S^{1,1}T^*M$ is the restriction of the natural projection  \eqref{eq:pi11} to $S^2T^*M$.
On the other hand, the corresponding almost Grassmannian structure $TM\cong E^*\otimes F$ yields 
\begin{equation*}
  S^2T^*M= (\La^2E\otimes \La^2F^*) \oplus (S^2E\otimes S^2F^*).
\end{equation*}
Analogously to Lemma \ref{l:decomp}, the two decompositions agree so that
$$
  \La^2E\otimes \La^2F^*\cong S^{1,1}(E^*\otimes F)
  \qtext{and}
  S^2E\otimes S^2F^*\cong\ker\pi^{1,1}.
$$
The metric on $M$ is {compatible} with the geometric structure if it is a section of $\La^2E\otimes \La^2F^*\cong S^{1,1}(E^*\otimes F)$.
In order that the metric is non-degenerate, the rank of the vector bundle $F$ has to be even.
Hence, if there is a compatible metric then the dimension of the base manifold is a multiple of 4.
It is also obvious, that all tangent vectors in the Segre cone are null with respect to any compatible metric.
Consequently, the compatible metric is of split signature.

If further the Levi-Civita connection of a compatible metric is a compatible connection of the geometric structure, then the metric (as well as the structure itself) is called \textit{para-quaternionic K\"ahler}.
Since Levi-Civita connection is torsion free, para-quaternionic K\"ahler structures are integrable.
It is also the case, that para-quaternionic K\"ahler metrics are necessarily Einstein.
In the 4-dimensional case, 
this feature may be stated so that the corresponding conformal manifold is anti-self-dual and contains an Einstein metric in the conformal class.
It follows that existence of para-quaternionic K\"ahler metrics is controlled by solutions to an invariant overdetermined system of differential equations, the so-called first BGG equation.
This is dealt in \cite{Bailey1991} in the holomorphic category with some minor additional assumptions.

\end{document}